\documentclass[12pt]{article}
\usepackage{amsmath}
\usepackage{graphicx,psfrag,epsf,color}
\usepackage{enumerate}
\usepackage{natbib}
\usepackage{url} 
\usepackage{amssymb}
\usepackage{multirow}
\usepackage{array}
\usepackage{bm}
\usepackage{booktabs}
\usepackage{float}
\usepackage[font=small,labelfont=bf]{caption}

\usepackage[total={6.4in,8.6in}]{geometry}
\newcommand{\blind}{0}

\newcommand{\Mean}{{\mbox{E}}}

\newcommand{\Cov}{{\mbox{cov}}}

\newcommand{\prob}{{\mbox{Pr}}}
\newtheorem{thm}{Theorem}[section]
\newtheorem{coro}{Corollary}[section]
\newtheorem{lemma}{Lemma}[section]
\newtheorem{defi}{Definition}[section]
\newtheorem{remark}{Remark}[section]

\begin{document}

\def\spacingset#1{\renewcommand{\baselinestretch}%
{#1}\small\normalsize} \spacingset{1}
%
%
\if0\blind
{
			\title{\bf A Massive Data Framework for M-Estimators with Cubic-Rate}
			\author{Chengchun Shi, Wenbin Lu and Rui Song\thanks{
					The authors gratefully acknowledge \textit{please remember to list all relevant funding sources in the unblinded version}}\hspace{.2cm}\\
				Department of Statistics, North Carolina State University 
			}
			\maketitle
	} \fi

\bigskip

\begin{abstract}
The divide and conquer method is a common strategy for handling massive data. In this article, we study the divide and conquer method for cubic-rate estimators under the massive data framework. We develop a general theory for establishing the asymptotic distribution of the aggregated M-estimators using a simple average. Under certain condition on the growing rate of the number of subgroups, the resulting aggregated estimators are shown to have faster convergence rate and asymptotic normal distribution, which are more tractable in both computation and inference than the original M-estimators based on pooled data. 
Our theory applies to a wide class of M-estimators with cube root convergence rate, including the location estimator, maximum score estimator and value search estimator. Empirical performance via simulations also validate our theoretical findings. 
\end{abstract}

\noindent%
\vfill

\section{Introduction}
\label{sec:intro}

In a world of explosively large data, effective estimation procedures are needed to deal with the computational challenge arisen from analysis of massive data. The divide and conquer method is a commonly used approach for handling massive data, which divides data into several groups and aggregate all subgroup estimators by a simple average to lessen the computational burden. A number of problems have been studied for the divide and conquer method, including variable selection \citep{Xie2014}, nonparametric regression \citep{zhang2013, zhao2016} and bootstrap inference \citep{Kleiner2014}, to mention a few. Most papers establish that the aggregated estimators achieve the oracle result, in the sense that they possess the same nonasymptotic error bounds or limiting distributions as the pooled estimators, which are obtained by fitting all the data in a single model. This implies that the divide and conquer scheme can not only maintain efficiency, but also obtain a feasible solution for analyzing massive data. 

In addition to the computational advantages for handling massive data, the divide and conquer method, somewhat surprisingly, can lead to aggregated estimators with improved efficiency over pooled estimators with slower than the usual $n^{1/2}$ convergence rate. There is a wide class of M-estimators with $n^{1/3}$ convergence rate. For example, \cite{Cher1964} studied a cubic-rate estimator for estimating the mode. It was shown therein that the estimator converges in distribution to the argmax of a Brownian motion minus a quadratic drift. \cite{Kim1990} systematically studied a class of cubic-rate M-estimators and established their limiting distributions as the argmax of a general Gaussian process minus a quadratic form. These results were extended to a more general class of M-estimators using modern empirical process results \citep{van1996,kosorok}. In this paper, we mainly focus on M-estimators with cubic-rate and develop a unified inference framework for the aggregated estimators obtained by the divide and conquer method. Our theory states that the aggregated estimators can achieve a faster convergence rate than the pooled estimators and have asymptotic normal distributions when the number of groups diverges at a proper rate as the sample size of each group grows. This enables a simple way for estimating the covariance matrix of the aggregated estimators.

When establishing the asymptotic properties of the aggregated estimators, a major technical challenge is to quantify the accumulated bias. Different from estimators with standard $n^{1/2}$ convergence rate, M-estimators with $n^{1/3}$ convergence rate generally do not have a nice linearization representation and the magnitude of the associated biases is difficult to quantify. In  literature, a few works have been developed for  studying the mean of the argmax of a simple one-dimensional Brownian motion process plus a deterministic function \citep[see for example,][]{Daniels1985, Cator2006, Pimentel2014}. In particular, \cite{Groene1999} 
provided a coupling inequality for the inverse process of the Grenander estimator based on Komlos-Major-Tusnady (KMT) approximation \citep{KMT1975}. However, it remains unclear and can be challenging to extend their technique under a more general setting. On one hand, the KMT approximation requires the underlying class of functions to be uniformly bounded \citep[see for example, ][]{Rio1994, KMT1994}. This assumption is violated in some applications for M-estimators, for example the value search estimator described in Section \ref{secexample}. On the other hand, their coupling inequality relies heavily on the properties of the argmax of a Brownian motion process with a parabolic drift \citep{Groene1989}, and is not applicable to cubic-rate estimators that converge to the argmax of a more general Gaussian process minus a quadratic form. Here, we propose a novel approach to derive a nonasymptotic error bound for the bias of aggregated M-estimators.

A key innovation in our analysis is to introduce a linear perturbation in the empirical objective function. 
In that way, we transform the problem of quantifying the bias into comparison of the expected supremum of the empirical objective function and that of its limiting Gaussian process. To bound the difference of these expected suprema,
we adopt similar techniques that have been recently studied by \cite{Cherno2013} and \cite{Cherno2014}. Specifically, they compared a function of the maximum for sum of mean-zero Gaussian random vectors with that of multivariate mean-zero random vectors with the same covariance function, and provided an associated coupling inequality. 
We improve their arguments by providing more accurate approximation results (Lemma \ref{lemmamaximasumsofgaussian}) for the identity function of maximums as needed in our applications. 

Another major contribution of this paper is to provide a tail inequality for cubic-rate M-estimators (Theorem \ref{thmconcenhath}). This helps us to construct a truncated estimator with bounded second moment, which is essential to apply Lindberg's central limit theorem for establishing the normality of the aggregated estimator. Under some additional tail assumptions on the underlying empirical process, our results can be viewed as a generalization of empirical process theories that establish consistency and $n^{1/3}$ convergence rate for the M-estimators.  Based on the results, we show that the asymptotic variance of the aggregated estimator can be consistently estimated by the sample variance of individual M-estimators in each group, which largely simplifies the inference procedure for M-estimators. 

The rest of the paper is organized as follows. We describe the divide and conquer method for M-estimators and state the major central limit theorem (Theorem \ref{thmcuberootdivideandconquer}) in Section \ref{sec:meth}. Three examples for the location estimator, maximum score estimator and value search estimator are presented in Section \ref{secexample} to illustrate the application of Theorem \ref{thmcuberootdivideandconquer}. Simulation studies are conducted in Section \ref{secsimu} to demonstrate the empirical performance of the aggregated estimators. Section \ref{secconceninequality} studies a tail inequality and Section \ref{secbias} provides the analysis of bias of M-estimators that are needed to prove Theorem \ref{thmcuberootdivideandconquer},
followed by a Discussion Section. All the technical proofs are provided in the Appendix.

\section{Method}\label{sec:meth}

The divide and conquer scheme for M-estimators is described as follows. In the first step, the data are randomly divided into several groups. For the $j$th group, consider the following M-estimator
\begin{eqnarray*}
	\hat{\theta}^{(j)}=\arg \max_{\theta \in \Theta} \mathbb{P}_{n_j}^{(j)} m(\cdot, \theta)\equiv\arg \max_{\theta \in \Theta}\frac{1}{n_j} \sum_{i=1}^{n_j} m(X_i^{(j)}, \theta), \qquad j=1,\dots,S,
\end{eqnarray*}
where $(X_1^{(j)},\dots, X_{n_j}^{(j)})$ denote the data for the $j$th group, $n_j$ is the number of observations in the $j$th group, $S$ is the number of groups, $m(\cdot, \cdot)$ is the objective function and $\theta$ is a $d$-dimensional vector of parameters that belong to a compact parameter space $\Theta$. In the second step, the aggregated estimator $\hat{\theta}_0$ is obtained as a simple average of all subgroup estimators, i.e.
\begin{eqnarray}\label{hattheta0}
	\hat{\theta}_0=\frac{1}{S}\sum_{j=1}^S \hat{\theta}^{(j)}.
\end{eqnarray}

We assume that all the $X_i^{(j)}$'s are independent and identically distributed across $i$ and $j$. In addition, for notational simplicity and without loss of generality, we assume equal allocation among $S$ groups, i.e. $n_1=\dots=n_S=n$. 
Here, we only consider M-estimation with non-smooth functions $m(\cdot, \theta)$ of $\theta$, and the resulting M-estimators $\hat{\theta}^{(j)}$'s have a convergence rate of $O_p(n^{-1/3})$. Such cubic-rate M-estimators have been widely studied in the literature, for example, the location estimator and maximum score estimator as demonstrated in the next section. The limiting distributions of these estimators have also been established using empirical process arguments \citep[cf.][]{Kim1990, van1996}. To be specific, let $\theta_0$ denote the unique maximizer of $\Mean \{m(\cdot, \theta)\}$, and assume $\theta_0 \in \Theta$. Then,
$\hat{h}^{(j)}\equiv n^{1/3} (\hat{\theta}^{(j)}-\theta_0)$ converges in distribution to $h_0 = \arg\max_h Z(h)$, where
\begin{eqnarray}\label{limitingprocess}
	Z(h)=G(h)-\frac{1}{2}h^T V h,
\end{eqnarray}
for some mean-zero Gaussian process $G$ and positive definite matrix $V=\partial^2 \Mean \{m(\cdot,\theta)\}/\partial \theta\theta^T|_{\theta=\theta_0}$.

The main goal of this paper is to establish the convergence rate and asymptotic normality of $\hat{\theta}_0$ under suitable conditions for $S$ and $n$, even though each $\hat{\theta}^{(j)}$ does not have a tractable limiting distribution. The dimension $d$ is assumed to be fixed, while the number of groups $S \rightarrow \infty$ as $n \rightarrow \infty$. Let $||\cdot||_2$ denote the Euclidean norm for vectors or induced matrix $L_2$ norm for matrices. We first introduce some conditions.

\medskip

\noindent (A1.) There exists a small neighborhood $N_\delta=\{\theta:||\theta-\theta_0||_2\le \delta\}$ 
in which $\Mean m\{(\cdot, \theta)\}$ is twice continuously differentiable with the Hessian matrix $-V(\theta)$, where $V(\theta)$ is positive definite in $N_\delta$. Moreover, assume $\Mean\{ m(\cdot,\theta_0)\}>\sup_{\theta \in N_\delta^c} \Mean \{m(\cdot, \theta)\}$.

\noindent (A2.) For any $\theta_1,\theta_2\in N_\delta$, we have $\Mean\{ |m(\cdot, \theta_1)-m(\cdot,\theta_2)|^2\}=O(||\theta_1-\theta_2||_2)$. 

\noindent (A3.) There exists some envelope function  $M(\cdot)\ge |m(\cdot, \theta)|$ for any $\theta$, and $\omega=||M(\cdot)||_{\psi_1}<\infty$, where $||\cdot||_{\psi_p}$ denotes the Orzlic norm of a random variable.

\noindent (A4.) The envelope function $M_R(\cdot)\equiv\sup_\theta \left\{\left|m(\cdot,\theta)\right|:||\theta-\theta_0||_2\le R\right\} $ satisfies $\Mean M_R^2=O(R)$ when $R\le \delta$. 

\noindent (A5.) The set of functions $\{m(\cdot, \theta)|\theta\in \Theta\}$ has Vapnik-Chervonenkis (VC) index $1\le v<\infty$.

\noindent (A6.) For any $\theta\in N_\delta$, $||V(\theta)-V||_2=O(||\theta-\theta_0||_2)$, where $V=V(\theta_0)$. 

\noindent (A7.) Let $L(\cdot)$ denote the variance process of $G(\cdot)$ satisfying $L(h)> 0$ whenever $h\neq 0$. (i) The function $L(\cdot)$ is symmetric and continuous, and has the rescaling property: $L(kh)=kL(h)$ for $k>0$. (ii) For any $h_1,h_2\in\mathbb{R}^d$ satisfying $||h_1||_2\le n^{1/3} \delta$ and $||h_2||_2\le n^{1/3} \delta$, we have
\begin{eqnarray*}
	\left|L(h_1-h_2)-n^{1/3} \Mean \left\{m(\cdot,\theta_0+n^{-1/3}h_1)-m(\cdot,\theta_0+n^{-1/3} h_2)\right\}^2\right|=O\left(\frac{(||h_1||+||h_2||)^2}{n^{1/3}}\right).
\end{eqnarray*}

\begin{thm}\label{thmcuberootdivideandconquer}
	Under Conditions (A1)-(A7), if $S=o(n^{1/6}/\log^{4/3} n)$ and $S\to \infty$ as $n\to \infty$, we have
	\begin{eqnarray}\label{hattheta0-theta0}
		\sqrt{S} n^{1/3} (\hat{\theta}_0-\theta_0) \stackrel{d}{\to} N(0, A),
	\end{eqnarray}
	for some positive definite matrix $A$. 
\end{thm}

\begin{remark}
	Theorem \ref{thmcuberootdivideandconquer} suggests that $\hat{\theta}_0$ converges at a rate of $O_p(S^{-1/2} n^{-1/3})$. In contrast, the original M-estimator obtained based on pooled data has a convergence rate of $O_p(S^{-1/3} n^{-1/3})$. This implies that we can gain efficiency by adopting the split and conquer scheme for cubic-rate M-estimators. Such result is interesting as most aggregated estimators in the divide and conquer literature share the same convergence rates as the original estimators based on pooled data. 
\end{remark}

\begin{remark}
	The constraints on $S$ suggest that the number of group cannot diverge too fast. A main reason as we showed in the proof of Theorem \ref{thmcuberootdivideandconquer} is that if $S$ grows too fast,  the asymptotic normality of $\hat{\theta}_0$ will fail due to accumulation of bias in the aggregation of subgroup estimators. 
	Given a data of size $N$, we can take $S\approx N^l$, $n=N/S\approx N^{1-l}$ with $l< 1/7$ to fulfill this requirement.
\end{remark}

\begin{remark}
	Conditions A1 - A5 and A7 (i) are similar to those in \cite{Kim1990} and are used to establish the cubic-rate convergence of the M-estimator in each group. Conditions A6 and A7 (ii) are used to establish the normality of the aggregated estimator. In particular, Condition A7 (ii) implies that 
	the Gaussian process $G(\cdot)$ has stationary increments, i.e. $\Mean[\{G(h_1) - G(h_2)\}^2] = L(h_1-h_2)$ for any $h_1,h_2\in\mathbb{R}^d$, which is used to control the bias of the aggregated estimator.
	
\end{remark}

In the rest of this section, we give a sketch for the proof of Theorem \ref{thmcuberootdivideandconquer}. The details of the proof are given in Sections \ref{secconceninequality} and \ref{secbias}. By the definitions of $\hat{\theta}_0$ and $\hat{h}^{(j)}$, it is equivalent to show
\begin{eqnarray}\label{clthattheta0}
	\frac{1}{\sqrt{S}} \sum_{j=1}^S \hat{h}^{(j)} \stackrel{d}{\to} N(0, A).
\end{eqnarray} 

When $S$ diverges, intuitively, \eqref{clthattheta0} follows by a direct application of Lindberg's central limit theorem for triangular arrays \citep[cf. Theorem 11.1.1,][]{Lahiri2006}. However, a few challenges remain. First, the estimator $\hat{h}^{(j)}$ may not possess finite second moment. Analogous to Kolmogorov's $3$-series theorem \citep[cf. Theorem 8.3.5,][]{Lahiri2006}, we handle this by first defining $\tilde{h}^{(j)}$, which is a truncated version of $\hat{h}^{(j)}$ with $||\tilde{h}^{(j)}||_2\le \delta_n$ 
for some $\delta_n>0$, such that $\sum_j \hat{h}^{(j)}$ and $\sum_j \tilde{h}^{(j)}$ are tail equivalent, i.e.
\begin{eqnarray*}
	\lim_k \prob\left(\displaystyle \bigcap_{n\ge k} \left\{\sum_{j=1}^{S(n)} \hat{h}^{(j)} = \sum_{j=1}^{S(n)}\tilde{h}^{(j)}\right\} \right)=1.
\end{eqnarray*}
Using Borel-Cantelli lemma, it suffices to show
\begin{eqnarray}\label{tildeh-hath}
	\sum_n \prob\left(\displaystyle \sum_{j=1}^{S(n)} \hat{h}^{(j)} = \sum_{j=1}^{S(n)}\tilde{h}^{(j)} \right)<\infty.
\end{eqnarray}

Now it remains to show
\begin{eqnarray*}
	\frac{1}{\sqrt{S}} \sum_{j=1}^S \tilde{h}^{(j)} = \frac{1}{\sqrt{S}} \sum_{j=1}^S \left\{\tilde{h}^{(j)} -\Mean (\tilde{h}^{(j)})\right\} + \sqrt{S}\Mean \tilde{h}^{(j)} \stackrel{d}{\to} N(0, A).
\end{eqnarray*}
The second challenge is to control the accumulated bias in the aggregated estimator, i.e.  showing
\begin{eqnarray}\label{bias}
	\sqrt{S}\Mean (\tilde{h}^{(j)})  \to 0.
\end{eqnarray}
Finally, it remains to show that the second moment of $\tilde{h}^{(j)}$ satisfies
\begin{eqnarray}\label{variance}
	\Mean (a^T \tilde{h}^{(j)})^2\to a^T A a,
\end{eqnarray}
for any $a\in \mathbb{R}^d$. Then, Theorem \ref{thmcuberootdivideandconquer} holds when \eqref{tildeh-hath}, \eqref{bias} and \eqref{variance} are established. Section \ref{secconceninequality} is devoted to verifying \eqref{tildeh-hath} and \eqref{variance}, while Section \ref{secbias} is devoted to proving \eqref{bias}. 

\section{Applications}\label{secexample}
In this section, we illustrate our main theorem (Theorem \ref{thmcuberootdivideandconquer}) with three applications including simple one-dimensional location estimator (Example \ref{locationest}) and more complicated multi-dimensional estimators with some constraints, such as maximum score estimator (Example \ref{scoreest}) and value-search estimator (Example \ref{valueest}).

\subsection{Location estimator}\label{locationest}
Let $X_i^{(j)}$ $(i=1,\dots,n; j=1,\dots,S)$ be i.i.d. random variables on the real line, with a continuous density $p$. In each subgroup $j$, consider the location estimator
\begin{eqnarray*}
	\hat{\theta}^{(j)}=\arg\max_{\theta \in \mathbb{R}} \frac{1}{n}\sum_{i=1}^n I(\theta-1\le X_i\le \theta+1).
\end{eqnarray*}

It was shown in Example 3.2.13 of \cite{van1996} and Example 6.1 of \cite{Kim1990} that each $\hat{\theta}^{(j)}$ has a cubic-rate convergence. 
We assume that $\prob(X\in [\theta-1,\theta+1])$ has a unique maximizer at $\theta_0$. When the derivative of $p$ exists and is continuous, $p'(\theta_0-1)-p'(\theta_0+1)>0$ implies that the second derivative of $\prob(X\in [\theta-1,\theta+1])$ is negative for all $\theta$ within some small neighbor $N_\delta$ around $\theta_0$. Therefore, Condition (A2) holds, since
\begin{eqnarray*}
	&&\Mean |I(\theta_1-1\le X\le \theta_1+1)-I(\theta_2-1\le X\le \theta_2+1)|^2\\
	&=&\prob(\theta_1-1\le X\le\theta_2-1)+\prob(\theta_1+1\le X\le \theta_2+1)\\
	&\le& \sup_{\theta\in N_\delta} \{p'(\theta-1)+p'(\theta+1)\} |\theta_1-\theta_2|,
\end{eqnarray*} 
for $\theta_1\le \theta_2$ and $|\theta_1-\theta_2|<0.5$. 
Moreover, if we further assume that $p$ has continuous second derivative in the neighborhood $N_\delta$, Condition (A6) is satisfied.    

The class of functions $\{|I(\theta-1\le X\le\theta+1)|: \theta \in \Theta\}$ is bounded by $1$ and belongs to VC class. In addition, we have
\begin{eqnarray*}
	&& \sup_{|\theta-\theta_0|<\epsilon} \left|I(\theta-1\le X\le\theta+1)-I(\theta_0-1\le X\le\theta_0+1)\right|\\
	&\le & I(\theta_0-1-\epsilon\le X\le \theta_0-1+\epsilon)+I(\theta_0+1-\epsilon\le X\le \theta_0+1+\epsilon),
\end{eqnarray*}
for small $\epsilon$. The $L_2(P)$ norm of the function on the second line is $O(\sqrt{\epsilon})$. Hence, Conditions (A4) and (A5) hold. 

Next, we claim that Condition (A7) holds for function $L(h)\equiv 2p(\theta_0+1)|h|$, or equivalently $\{p(\theta_0-1)+p(\theta_0+1)\}|h|$, since $p(\theta_0-1)=p(\theta_0+1)$. Obviously, $L(\cdot)$ is symmetric and satisfies the rescaling property. For any $h_1,h_2$ such that $\max(|h_1|,|h_2|)\le n^{1/3}\delta$, we define $\theta_1=\theta_0+n^{-1/3}h_1 \in N_\delta$ and $\theta_2=\theta_0+n^{-1/3}h_2\in N_\delta$. Let $[a,b]$ denote the indicator function $I(a\le X\le b)$. Assume $h_1\le h_2$. We have
\begin{eqnarray*}
	&&n^{1/3}\Mean \left|[\theta_1-1,\theta_1+1]-[\theta_2-1,\theta_2+1]\right|^2=n^{1/3}\Mean [\theta_1-1,\theta_2-1]+n^{1/3}\Mean[\theta_1+1,\theta_2+1]\\
	&=&n^{1/3}\int_{\theta_1-1}^{\theta_2-1} p(\theta)d\theta+n^{1/3}\int_{\theta_1+1}^{\theta_2+1}p(\theta)d\theta=\{p(\theta_0+1)+p(\theta_0-1)\}(h_2-h_1)+R,
\end{eqnarray*}
where the remainder term $R$ is bounded by
\begin{eqnarray*}
	&&\sup_{\theta_1\le \theta \le \theta_2}\left(|p(\theta-1)-p(\theta_0-1)|+|p(\theta+1)-p(\theta_0+1)|\right) (h_2-h_1)\\
	&\le & \sup_{\theta \in N_\delta} 4n^{-1/3}|p'(\theta)| (h_2-h_1)  \max(|h_1|, |h_2|)\le \sup_{\theta \in N_\delta} 4n^{-1/3}|p'(\theta)| (|h_1|+|h_2|)^2,
\end{eqnarray*}
using a first order Taylor expansion. The case when $h_1>h_2$ can be similarly discussed. Therefore, Condition (A7) holds. Theorem \ref{thmcuberootdivideandconquer} then follows.

\subsection{Maximum score estimator}\label{scoreest}
Consider the regression model $Y_i^{(j)}={X_i^{(j)}}^T \beta_0+e_i^{(j)}$, $,j=1,\cdots,S$, where $X_i^{(j)}$ is a $d$-dimensional vector of covariates and $e_i^{(j)}$ is the random error. Assume that $(X_i^{(j)}, e_i^{(j)})$'s are i.i.d. copies of $(X, e)$. The maximum score estimator is defined as 
\begin{eqnarray*}
	\hat{\beta}^{(j)}=\arg\max_{||\beta||_2= 1} \sum_{i=1}^n \{I(Y_i^{(j)}\ge 0, {X_i^{(j)}}^T \beta\ge 0)+I(Y_i^{(j)}< 0, {X_i^{(j)}}^T \beta< 0)\},
\end{eqnarray*}
where the constraint $||\beta||_2=1$ is to guarantee the uniqueness of the maximizer. 

Assume $||\beta_0||=1$, otherwise we can define $\beta^\star=\beta_0/||\beta_0||_2$ and establish the asymptotic distribution of $\hat{\beta}_0-\beta^\star$ instead. It was shown in Example 6.4 of \cite{Kim1990} that $\hat{\beta}^{(j)}$ has a cubic-rate convergence, when (i) $\hbox{median}(e|X)=0$;
(ii) $X$ has a bounded, continuously differentiable density $p$; and (iii) the angular component 
of $X$ has a bounded continuous density with respect to the surface measure on $\mathbb{S}^{d-1}$, which corresponds to the unit sphere in $\mathbb{R}^d$. 

Theorem \ref{thmcuberootdivideandconquer} is not directly applicable to this example since Assumption (A1) is violated. The Hessian matrix
\begin{eqnarray*}
V=-\frac{\partial^2 \Mean \{I(Y_i^{(j)}\ge 0, {X_i^{(j)}}^T \beta\ge 0)+I(Y_i^{(j)}< 0, {X_i^{(j)}}^T \beta< 0)\}}{\partial \beta \beta^T}\rvert_{\beta_0}
\end{eqnarray*}	
is not positive definite.
One possible solution is to use the arguments from the constrained M-estimator literature \citep[e.g.][]{Geyer1994} to approximate the set $||\beta||_2=1$ by a hyperplane $(\beta-\beta_0)^T \beta=0$, and obtain a version of Theorem \ref{thmcuberootdivideandconquer} for the constrained cubic-rate M-estimators. We adopt an alternative approach here, and consider a simple reparameterization to make Theorem \ref{thmcuberootdivideandconquer} applicable.


By Gram-Schmidt orthogonalization, we can obtain an orthogonal matrix $[\beta_0, U_0]$ with $U_0$ being a $\mathbb{R}^{d\times (d-1)}$ matrix subject to the constraint $U_0^T \beta_0=0$. Define
\begin{eqnarray}\label{betathetafun}
	\beta(\theta)=\sqrt{1-||\theta||_2^2} \beta_0+U_0\theta,
\end{eqnarray}
for all $\theta\in \mathbb{R}^{d-1}$ and $||\theta||_2\le 1$. Take $\Theta$ to be the unit ball $B_2^{d-1}$ in $\mathbb{R}^{d-1}$. Define
\begin{eqnarray*}
	\hat{\theta}^{(j)}=\arg\max_{\theta \in \Theta} \sum_{i=1}^n [I(Y_i^{(j)}\ge 0, {X_i^{(j)}}^T \beta(\theta) \ge 0)+I(Y_i^{(j)}< 0, {X_i^{(j)}}^T \beta(\theta)< 0)].
\end{eqnarray*}
Note that under the assumption $\hbox{median}(e|X)=0$, we have $\theta_0=0$.

Let $m(y, x, \beta)=I(y\ge 0, x^T \beta\ge 0)+I(y< 0, x^T \beta< 0)$. 
Define
\begin{eqnarray*}
	\kappa(x)=\Mean \{I(e+X^T \beta_0\ge 0)-I(e+X^T\beta_0<0) |X=x\}.
\end{eqnarray*}
It is shown in \cite{Kim1990} that
\begin{eqnarray}\label{firstderivative}
	\frac{\partial \Mean \{m(\cdot,\cdot,\beta)\}}{\partial \beta}=||\beta||_2^{-2} \beta^T \beta_0 (I+||\beta||_2^{-2} \beta \beta^T) \int_{x^T \beta_0=0} \kappa(T_\beta x)p(T_\beta x)d\sigma,
\end{eqnarray}
where
\begin{eqnarray*}
	T_\beta=(I-||\beta||_2^{-2} \beta \beta^T) (I-\beta_0 \beta_0^T) +||\beta||_2^{-1} \beta \beta_0^T,
\end{eqnarray*}
and $\sigma$ is the surface measure on the line $x^T \beta_0=0$. 

Note that $\partial \beta(\theta)/\partial \theta$ has finite derivatives for all orders as long as $||\theta||_2 < 1$. Assume that $\kappa$ and $p$ have twice continuous derivatives. This together with \eqref{firstderivative} implies that $\Mean\{ m(\cdot,\cdot,\beta(\theta))\}$ has third continuous derivative as a function of $\theta$ in a small neighborhood $N_\delta$ ($\delta<1$) around $0$. This verifies (A6). Moreover, for any $\theta_1$, $\theta_2\in N_\delta$ with $||\theta_1-\theta_2||_2\le \epsilon$, we have 
\begin{eqnarray}\nonumber
	&&||\beta(\theta_1)-\beta(\theta_2)||_2^2=||\theta_1-\theta_2||_2^2+\left(\sqrt{1-||\theta_1||_2^2}-\sqrt{1-||\theta_2||_2^2}\right)^2\\ \label{betathetaderivative}
	&=&||\theta_1-\theta_2||_2^2+\frac{\left(1-||\theta_1||_2^2-1+||\theta_2||_2^2\right)^2}{\left(\sqrt{1-||\theta_1||_2^2}+\sqrt{1-||\theta_2||_2^2}\right)^2}\le \frac{2||\theta_1-\theta_2||_2^2}{1-\delta^2}. 
\end{eqnarray} 
\cite{Kim1990} showed that $\Mean\{ |m(\cdot,\cdot,\beta_1)-m(\cdot,\cdot,\beta_2)|\}=O(||\beta_1-\beta_2||_2)$ near $\beta_0$. This together with \eqref{betathetaderivative} implies
\begin{eqnarray*}
	\Mean\{ |m(\cdot,\cdot,\beta(\theta_1))-m(\cdot,\cdot,\beta(\theta_2))|^2\}\le 2\Mean\{ |m(\cdot,\cdot,\beta(\theta_1))-m(\cdot,\cdot,\beta(\theta_2))|\}=O(||\theta_1-\theta_2||_2).
\end{eqnarray*}
Therefore, (A2) is satisfied and (A3) trivially holds since $|m|\le 1$. 

It was also shown in \cite{Kim1990} that the envelope $M_\epsilon$ of the class of functions $\{m(\cdot,\cdot,\beta)-m(\cdot,\cdot,\beta_0): ||\beta-\beta_0||_2\le \epsilon\}$ satisfies $\Mean M_\epsilon^2=O(\epsilon)$. Using \eqref{betathetaderivative}, we can show that the envelope $\tilde{M}_\epsilon$ of the class of functions $\{m(\cdot,\cdot,\beta(\theta))-m(\cdot,\cdot,\beta_0): ||\theta||_2\le \epsilon\}$ also satisfies $\Mean \tilde{M}_\epsilon^2=O(\epsilon)$. Thus, (A4) is satisfied. Moreover, since the class of functions $m(\cdot,\cdot,\beta)$ over all $\beta$ belongs to the VC class, so does the class of function $m(\cdot, \cdot, \beta(\theta))$. This verifies (A5). 

Finally, we establish (A7). For any $\theta_1,\theta_2\in N_\delta$, define $h_1=n^{1/3} \theta_1$ and $h_2=n^{1/3} \theta_2$. We have
\begin{eqnarray}\nonumber
	&&n^{1/3}\Mean\left\{ \left|m(Y,X,\beta(h_1/n^{1/3}))-m(Y,X,\beta(h_2/n^{1/3}))\right|^2\right\}\\ \nonumber
	&=&n^{1/3} \Mean \left\{\left|I(X^T \beta(h_1/n^{1/3})\ge 0)-I(X^T \beta(h_2/n^{1/3})\ge 0)\right|I(Y\ge 0)\right\}\\ \nonumber
	&+& n^{1/3}\Mean \left\{\left|I(X^T \beta(h_1/n^{1/3})< 0)-I(X^T \beta(h_2/n^{1/3})<0)\right|I(Y<0)\right\}\\ \label{E2A7eq1}
	&=&n^{1/3} \Mean\left\{ \left|I(X^T \beta(h_1/n^{1/3})\ge 0)-I(X^T \beta(h_2/n^{1/3})\ge 0)\right|\right\}.
\end{eqnarray}
We write $X$ as $r\beta_0+z$ with $z$ orthogonal to $\beta_0$. Equation \eqref{E2A7eq1} can be written as
\begin{eqnarray}\label{E2A7eq2}
	n^{1/3} \Mean \left\{ \left|I\left(r\sqrt{1-\left|\left|\frac{h_1}{n^{1/3}}\right|\right|_2^2}+z^T U \frac{h_1}{n^{1/3}}\ge 0\right)-I\left(r\sqrt{1-\left|\left|\frac{h_2}{n^{1/3}}\right|\right|_2^2}+z^T U \frac{h_2}{n^{1/3}}\ge 0\right)\right|\right\}.
\end{eqnarray}
Define $\omega=n^{1/3} r$. Equation \eqref{E2A7eq2} can be expressed as
\begin{eqnarray*}
	\int \int I(-z^T U h_1 (1-n^{-2/3} ||h_1||_2^2)^{-1/2}>\omega \ge -z^T U h_2 (1-n^{-2/3} ||h_2||_2^2)^{-1/2})p\left(\frac{\omega}{n^{1/3}}, z\right)d\omega dz.
\end{eqnarray*}
Assume that $p(r, z)$ is differentiable with respect to $r$ and $|\partial p(r, z)/\partial r|\le q(z)$ for some function $q$. Then, \eqref{E2A7eq2} is equal to
\begin{eqnarray*}
	&&\int |z^T U \{h_1 (1-n^{-2/3} ||h_1||_2^2)^{-1/2}-h_2 (1-n^{-2/3} ||h_2||_2^2)^{-1/2}\}| p(0, z)dz+R_1\\
	&=& \int |z^T U(h_1-h_2)|p(0, z) dz + R_1+R_2,
\end{eqnarray*}
where the remainders $|R_1|$ and $|R_2|$ are bounded by 
\begin{eqnarray*}
	|R_1|\le \int n^{-1/3} \{(z^T U h_1)^2 + (z^T U h_2)^2\} q(z) dz=O(n^{-1/3} \{||h_1||_2^2+||h_2||_2^2\}),
\end{eqnarray*}
and 
\begin{eqnarray*}
	|R_2|&\le& |(1-n^{-2/3} ||h_1||_2^2)^{-1/2}-1| \int |z^T U h_1| p(0, z)dz\\
	&+&|(1-n^{-2/3} ||h_2||_2^2)^{-1/2}-1| \int |z^T U h_2| p(0, z)dz\\
	&\le & n^{-1/3} (||h_1||_2+||h_2||) \int (|z^T U h_1|+ |z^T U h_2|)p(0,z)dz = O(n^{-1/3} \{||h_1||_2^2+||h_2||_2^2\}),
\end{eqnarray*}
under suitable moment assumptions on functions $p(0, z)$ and $q(z)$. This verifies (A7). 

An application of Theorem \ref{thmcuberootdivideandconquer} implies 
\begin{eqnarray*}
	\frac{1}{\sqrt{S}} \sum_{j=1}^S n^{1/3}\hat{\theta}^{(j)}\stackrel{d}{\to} N(0, A),
\end{eqnarray*}
for some positive definite matrix $A \in \mathbb{R}^{(d-1)\times (d-1)}$. Hence
\begin{eqnarray}\label{anormalUtheta}
	\frac{1}{\sqrt{S}} \sum_{j=1}^S n^{1/3} U \hat{\theta}^{(j)}\stackrel{d}{\to} N(0, U A U^T). 
\end{eqnarray}
By the definition of $\hat{\theta}^{(j)}$ and $\hat{\beta}^{(j)}$, we have
\begin{eqnarray*}
	&& \left|\frac{1}{\sqrt{S}} \sum_{j=1}^S n^{1/3} (\hat{\beta}^{(j)}-\beta_0-U \hat{\theta}^{(j)})\right|\le \left|\frac{1}{\sqrt{S}}\sum_{j=1}^S n^{1/3} |\sqrt{1-||\hat{\theta}^{(j)}||_2^2}-1|\right|\\
	&\le &\left|\frac{1}{\sqrt{S}}\sum_{j=1}^S n^{1/3} \frac{|1-||\hat{\theta}^{(j)}||_2^2-1|}{|\sqrt{1-||\hat{\theta}^{(j)}||_2^2}+1|}\right|\le \frac{n^{1/3}}{\sqrt{S}} \sum_j ||\hat{\theta}^{(j)}||_2^2.
\end{eqnarray*}
With probability at least $1-S/n\to 1$, the last expression is equal to $O(\sqrt{S} n^{1/3} n^{-2/3} \log^{2/3} n)=o(1)$, which is implied by the tail inequality for $\hat{\theta}^{(j)}$ established in Theorem \ref{thmconcenhath}. Combining this together with \eqref{anormalUtheta}, we have
\begin{eqnarray*}
	\frac{1}{\sqrt{S}} \sum_{j=1}^S n^{1/3} (\hat{\beta}^{(j)}-\beta_0)\stackrel{d}{\to} N(0, U A U^T). 
\end{eqnarray*} 

\subsection{Value search estimator}\label{valueest}
The value search estimator was introduced by \cite{zhang2012} for estimating the optimal treatment regime. The data can be summarized as i.i.d. triples $\{O_{i}^{(j)}=(X_{i}^{(j)}, A_{i}^{(j)}, Y_{i}^{(j)}), i=1,\dots,n; j=1,\dots,S\}$, where $X_{i}^{(j)}\in \mathbb{R}^d$ denote patient's baseline covariates, $A_{i}^{(j)}$ is the treatment received by the patient taking the value $0$ or $1$, and $Y_{i}^{(j)}$ is the response, the larger the better by convention. Consider the following model 
\begin{eqnarray}\label{model}
	Y_{i}^{(j)}=\mu(X_{i}^{(j)})+A_{i}^{(j)} C(X_{i}^{(j)})+e_{i}^{(j)},
\end{eqnarray}
where $\mu(\cdot)$ is the baseline mean function, $C(\cdot)$ is the contrast function, and $e_{i}^{(j)}$ is the random error with $E\{e_{i}^{(j)}|A_i^{(j)},X_i^{(j)}\} = 0$. The optimal treatment regime is defined in the potential outcome framework. Specifically, let $Y_{i}^{(j)\star}(0)$ and $Y_{i}^{(j)\star}(1)$ be the potential outcomes that would be observed if the patient received treatment $0$ or $1$, accordingly. For a treatment regime $d$ that maps $X_{i}^{(j)}$ to $\{0, 1\}$, define the potential outcome
\begin{eqnarray*}
	Y_{i}^{(j) \star}(d)=d(X_{i}^{(j)}) Y_{i}^{(j)\star}(1)+\{1-d(X_{i}^{(j)})\} Y_{i}^{(j)\star}(0).
\end{eqnarray*}
The optimal regime $d^{\tiny{opt}}$ is defined as the rule that maximizes the expected potential outcome, i.e, the value function,  $\Mean \{Y_{i}^{(j)\star}(d)\}$. Under the
stable unit treatment value assumption (SUTVA) and no unmeasured confounders assumption \citep{Neyman1990}, the optimal treatment regime under model \eqref{model} is given by $d^{\tiny{opt}}(x)=I\{C(x)>0\}$.

The true contrast function $C(\cdot)$ can be complex. As suggested by \cite{zhang2012}, in practice we can find the restricted optimal regimen within a class of decision rules, such as linear treatment decision rules 
$d(x,\beta)=I(\beta_1+x_1 \beta_2+\cdots+x_d \beta_{d+1}>0)$ indexed by $\beta \in \mathbb{R}^{d+1}$,
where the subscript $k$ denotes the $k$th element in the vector. Let $\beta^\star = \arg\max_{\beta}V(\beta)$, where 
$V(\beta)=\Mean \{Y_{i}^{(j)\star}(d(X_{i}^{(j)}, \beta))\}$.
To make $\beta^\star$ identifiable, 
we assume $\beta^\star_1=-1$. 
Define $\theta^\star = (\beta^\star_2,\cdots,\beta^\star_{d+1})^T$. The restricted optimal treatment regime is given by $\tilde{d}^{opt}(x,\theta^\star) = I(x^T\theta^\star > 1)$ and the value function is defined by $V(\theta)=\Mean \{Y_{i}^{(j)\star}(\tilde{d}(X_{i}^{(j)}, \theta))\}$ with $\tilde{d}(x, \theta)=I(x^T\theta>1)$. \cite{zhang2012} proposed an inverse propensity score weighted estimator of the value function $V(\theta)$ and the associated value search estimator by maximizing the estimated value function. Specifically, for each group $j$, the value search estimator is defined as
\begin{eqnarray*}
	\hat{\theta}^{(j)}=\arg\max_{\theta \in \Theta} \frac{1}{n}\sum_{i=1}^n \frac{\tilde{d}(X_{i}^{(j)}, \theta)A_{i}^{(j)}+\{1-\tilde{d}(X_{i}^{(j)}, \theta)\}(1-A_{i}^{(j)})}{\pi_{i}^{(j)}A_{i}^{(j)}+(1-\pi_i^{(j)})(1-A_{i}^{(j)})}Y_{i}^{(j)},
\end{eqnarray*}
where $\pi_i^{(j)}=\prob(A_i^{(j)}=1|X_i^{(j)})$ is the propensity score and known in a randomized study. Here, for illustration purpose, we assume that $\pi_i^{(j)}$'s are known. 

Define $m(O_i^{(j)}, \theta)=\xi_i^{(j)} \tilde{d}(X_i^{(j)}, \theta)$, where 
\begin{eqnarray*}
	\xi_i^{(j)}=\frac{A_i^{(j)}}{\pi_i^{(j)}} C(X_i^{(j)})+\frac{A_i^{(j)}-\pi_i^{(j)}}{\pi_i^{(j)} (1-\pi_i^{(j)})}\left\{\mu(X_i^{(j)})+e_i^{(j)}\right\}.
\end{eqnarray*}
With some algebra, we can show that $\hat{\theta}^{(j)}$ also maximizes $\mathbb{P}_n^{(j)} m(\cdot, \theta)$, where $\mathbb{P}_n^{(j)}$ is the empirical measure for data in group $j$. 
Unlike the previous two examples, here the function $m$ is not bounded. To fulfill (A3), we need $||\xi_i^{(j)}||_{\psi_1}<\infty$. This holds when $0<\gamma_1<\pi_i^{(j)}<\gamma_2<1$ for some constants $\gamma_1$ and $\gamma_2$, $||C(X_i^{(j)})||_{\psi_1}<\infty$, $||\mu(X_i^{(j)})||_{\psi_1}<\infty$ and $||e_i^{(j)}||_{\psi_1}<\infty$.


To show (A1) and (A6), we evaluate the integral
\begin{eqnarray}\label{someintegral}
	\Gamma(\theta)=\Mean \{\xi \tilde{d}(X, \theta)\}=\Mean \{C(X) \tilde{d}(X, \theta)\}=\int_{x^T \theta>1} C(x)p(x)dx,
\end{eqnarray}
where $p(x)$ is the density function of $X_i^{(j)}$. Consider the transformation 
\begin{eqnarray*}
	T_\theta=(I-||\theta||_2^{-2} \theta \theta^T)+||\theta||_2^{-2} \theta (\theta^\star)^T,
\end{eqnarray*}
which maps the region $\{x^T \theta^\star >1\}$ onto $\{x^T \theta>1 \}$, and $\{x^T \theta^\star=1 \}$ onto $\{x^T \theta=1\}$. We exclude the trivial case with $\theta^\star=0$. The above definition is meaningful when $\theta$ is taken over a small neighborhood $N_\delta$ of $\theta^\star$. We assume that functions $p$ and $C$ are continuously differentiable. Note that
\begin{eqnarray*}
	\frac{\partial T_\theta x}{\partial \theta}=-\frac{\{\theta^T x-(\theta^\star)^T x\}}{||\theta||_2^2}I-\frac{\theta x^T}{||\theta||_2^2}+\frac{2\theta \theta^T (x^T \theta-x^T \theta^\star)}{||\theta||_2^4}.
\end{eqnarray*}
Using some differential geometry arguments similarly as in Section 5 of \cite{Kim1990}, we can show that the integral \eqref{someintegral} can be represented as
\begin{eqnarray*}
	\Gamma(\theta)=\int_{x^T \theta^\star > 1} \left[-\frac{1}{||\theta||_2^2}\theta^T\frac{\partial C(x)p(x)}{\partial x}x+\frac{\{\theta^T x-(\theta^\star)^T x\} }{||\theta||_2^4} \theta^T \frac{\partial C(x)}{\partial x}\theta-\frac{\theta^T x-(\theta^\star)^T x}{||\theta||_2^2} \frac{\partial C(x) p(x)}{\partial x}\right]dx,
\end{eqnarray*}
which is thrice differentiable under certain conditions on $C(x)$, $p(x)$ and their derivatives. 

To show (A7), we assume that the conditional density $p(x|y)$ of $X$ given $Y=1-X^T\theta^\star$ exists and is continuously differentiable with respect to $y$. Similarly assume that the density $q(y)$ of $Y$  exists and is continuously differentiable. Let $g(X)=\Mean (\xi^2|X)$.  For any $h_1,h_2 \in \mathbb{R}^d$, we have
\begin{eqnarray*}
	&&n^{1/3} \Mean\left\{ \xi^2 \left|I(X^T \theta^\star+n^{-1/3} X^T h_1>1)-I(X^T \theta^\star+n^{-1/3} X^T h_2>1)\right|^2\right\}\\
	&=&n^{1/3} \int g(x) \left|I(n^{-1/3} x^T h_1>y)-I(n^{-1/3} x^T h_2>y)\right|p(x|y)q(y)dx dy.
\end{eqnarray*} 

Let $y=n^{-1/3} z$. The last expression in the above equation can be written as
\begin{eqnarray*}
	&&\int g(x) \left|I(x^T h_1>z)-I(x^T h_2>z)\right|p(x|0)q(0)dx dz+R\\
	&=& \int g(x) |x^T (h_1-h_2)| p(x|0) q(0) dx+R,
\end{eqnarray*}
with the remainder term 
\begin{eqnarray*}
	R=\int g(x) \left|I(x^T h_1>z)-I(x^T h_2>z)\right|\{p(x|n^{-1/3}z)q(n^{-1/3} z)-p(x|0)q(0)\}dx dz,
\end{eqnarray*}
which is $O(n^{-1/3} (||h_1||_2^2+||h_2||_2^2))$ under certain conditions on $q(x)$ and $p(x|\cdot)$. Conditions (A2) and (A4) can be similarly verified. Since the class of functions $\{g(x)I(x^T \theta>1):\theta \in \mathbb{R}^d \}$ has finite VC index, Condition (A5) also holds. Theorem \ref{thmcuberootdivideandconquer} then follows.

\section{Numerical studies}\label{secsimu}
In this section, we examine the numerical performance of the aggregated M-estimator for the three examples studied in the previous section and compare it with the M-estimator based on pooled data, denoted as the pooled estimator. 

\subsection{Location estimator}
The data $X_j$ $(j=1,\dots,N)$ were independently generated from the standard normal distribution. The true parameter $\theta_0$ that maximizes $\Mean \{I(\theta-1\le X_j\le \theta+1)\}$ was set to be $0$. Let $\tilde{\theta}_0$ and $\hat{\theta}_0$ denote the pooled estimator and the aggregated estimator, respectively. To obtain $\hat{\theta}_0$, we randomly divided the data into $S$ blocks with equal size $n=N/S$.

We took $N=2^i$ for $i=12$, 14, 16, 18, 20, and $S=2^j$ for $j=3,\dots,7$. For each combination of $N$ and $S$, we estimated the standard error of $\hat{\theta}_0$ by
$$
\widehat{SE}(\hat{\theta}_0)= \frac{1}{\sqrt{S}} \left\{\frac{1}{S-1}\sum_{l=1}^S \left(\hat{\theta}^{(l)}-\hat{\theta}_0\right)^2\right\}^{1/2},
$$
where $\hat{\theta}^{(l)}$ denotes the M-estimator for the $l$th group. For each scenario, we conducted $1000$ simulation replications and report the bias and sample standard deviation (denoted as SD) of estimators $\tilde{\theta}_0$ and $\hat{\theta}_0$, and estimated standard error and coverage probability (denoted as CP) of Wald-type 95\% confidence interval for $\hat{\theta}_0$ in Table \ref{tableE1}.  

Based on the results, it can be seen that the aggregated estimator $\hat{\theta}_0$ has much smaller standard deviation than the pooled estimator $\tilde{\theta}_0$, indicating the efficiency gain by the divide and conquer scheme as shown in our theory.  In addition, the bias of $\hat{\theta}_0$ generally becomes bigger and the standard deviation of $\hat{\theta}_0$ generally becomes smaller when $S$ and $N$ increase, and the normal approximation becomes more accurate when $S$ increases. This demonstrates the bias-variance trade off for aggregated estimators. With properly chosen $S$, the estimated standard error of $\hat{\theta}_0$ is close to its standard deviation and the coverage probability is close to the nominal level.

\begin{center}
	\captionof{table}{Bias and standard deviation (SD) for $\tilde{\theta}_0$ and $\hat{\theta}_0$, and estimated standard error and coverage probabilities (CP) of Wald-type 95\% confidence intervals for $\hat{\theta}_0$.}\label{tableE1} 
	\small
	\begin{tabular}{cccccccc}
		\toprule 
		& & $\tilde{\theta}_0$ & \multicolumn{5}{c}{$\hat{\theta}_0$}\\ 
		&&  & $S=2^3$ & $S=2^4$ & $S=2^5$ & $S=2^6$ & $S=2^7$\\ \midrule
		\multirow{4}{*}{$N=2^{12}$} & Bias & $-0.0017$ & $-0.0022$ & $-0.0036$ & $-0.0070$ & $-0.0162$ & $-0.0116$ \\ 
		& SD & $0.0549$ & $0.0397$ & $0.0350$ & $0.0339$ & $0.0297$ & $0.0496$ \\
		& $\widehat{SE}$ & NA & $0.0377$ & $0.0349$ & $0.0314$ & $0.0288$ & $0.0415$ \\ 
		& CP & NA  &$0.913$ & $0.932$& $0.919$ & $0.901$ & $0.831$ \\ \midrule 
		\multirow{4}{*}{$N=2^{14}$} & Bias & $-0.0032$ & $-0.0016$ & $0.0003$ & $-0.0015$ & $-0.0045$ & $-0.0073$ \\ 
		& SD & $0.0359$ & $0.0248$ & $0.0222$ & $0.0196$ & $0.0180$ & $0.0157$ \\
		& $\widehat{SE}$ & NA & $0.0239$ & $0.0218$ & $0.0197$ & $0.0175$ & $0.0159$ \\ 
		& CP & NA & $0.912$ & $0.922$ & $0.953$ & $0.940$ & $0.924$\\ \midrule 
		\multirow{4}{*}{$N=2^{16}$} & Bias & $-0.0008$ & $-0.0002$ & $-0.0004$ & $-3.7\times 10^{-5}$ & $-0.0010$ & $-0.0012$ \\ 
		& SD & $0.0216$ & $0.0164$ & $0.0141$ & $0.0129$ & $0.0111$ & $0.0097$ \\
		& $\widehat{SE}$ & NA & $0.0151$ & $0.0138$ & $0.0123$ & $0.0111$ & $0.0099$ \\
		& CP & NA & $0.900$ & $0.939$ & $0.932$ & $0.939$ & $0.946$\\ \midrule 
		\multirow{4}{*}{$N=2^{18}$} & Bias & $-0.0003$ & $-0.0005$ & $-0.0006$ & $-0.0004$ & $-0.0002$ & $7.5\times 10^{-5}$ \\ 
		& SD & $0.0141$ & $0.0101$ & $0.0089$ & $0.0078$ & $0.0070$ & $0.0064$ \\
		& $\widehat{SE}$ & NA & $0.0094$ & $0.0086$ & $0.0078$ & $0.0070$ & $0.0062$ \\
		& CP & NA & $0.899$ & $0.924$ & $0.940$ & $0.945$ & $0.943$\\ \midrule 
		\multirow{4}{*}{$N=2^{20}$} & Bias & $4.8\times 10^{-5}$ & $-9.7\times 10^{-5}$ & $6.0\times10^{-5}$ & $0.0003$ & $-0.0001$ & $-4.3\times 10^{-5}$ \\ 
		& SD & $0.0087$ & $0.0062$ & $0.0055$ & $0.0049$ & $0.0046$ & $0.0041$ \\
		& $\widehat{SE}$ & NA & $0.0060$ & $0.0055$ & $0.0049$ & $0.0044$ & $0.0040$ \\ 
		& CP & NA & $0.920$ & $0.940$ & $0.942$ & $0.935$ & $0.936$\\ \bottomrule
	\end{tabular}
\end{center}



\subsection{Maximum score estimator}

Consider the model $Y_i=1.5 X_{i1}-1.5 X_{i2}+0.5 e_i$, $i=1,\cdots,N$, where $X_{i1}$, $X_{i2}$ and $e_i$ were generated independently from the standard normal. Hence, $\theta_0=(\theta_1,\theta_2)^T=(1.5,-1.5)^T$. Let $\tilde{\theta}_0 = (\tilde{\theta}_1,\tilde{\theta}_2)^T$ denote the pooled estimator and $\hat{\theta}_0=(\hat{\theta}_1, \hat{\theta}_2)^T$ the aggregated estimator. We set $N=2^{18}, 2^{20}, 2^{22}$ and $S=2^4, \dots, 2^{7}$. Table \ref{tableE2} reports the results based on 1000 replications. The findings are very similar to those for the location estimator in the previous example.


\begin{center}
	\captionof{table}{Bias, standard deviation (SD), standard error ($\widehat{SE}$) for $\tilde{\theta}_1,\tilde{\theta}_2$ and $\hat{\theta}_1, \hat{\theta}_2$, and coverage probabilities of CI for $\hat{\theta}_1, \hat{\theta}_2$.}\label{tableE2}  
	\small
	\begin{tabular}{ccccccc}
		\toprule 
				& & $\tilde{\theta}_1$ & \multicolumn{4}{c}{$\hat{\theta}_1$}\\ 
		&&  & $S=2^4$ & $S=2^5$ & $S=2^6$ & $S=2^7$\\ \midrule
		\multirow{4}{*}{$N=2^{18}$} & Bias & $-1.3\times 10^{-5}$ & $4.6 \times 10^{-4}$ & $0.0011$ & $0.0023$ & $0.0040$ \\ 
		& SD & $0.0049$ & $0.0030$ & $0.0027$ & $0.0025$ & $0.0020$ \\
		& $\widehat{SE}$ & NA & $0.0030$ & $0.0027$ & $0.0024$ & $0.0020$ \\ 
		& CP & NA  &$0.932$ & $0.922$& $0.812$ & $0.467$ \\ \midrule 
		\multirow{4}{*}{$N=2^{20}$} & Bias & $-8.1\times 10^{-6}$ & $7.7\times 10^{-5}$ & $1.9\times 10^{-4}$ & $4.4\times 10^{-4}$ & $0.0013$ \\ 
		& SD & $0.0032$ & $0.0020$ & $0.0017$ & $0.0014$ & $0.0014$ \\
		& $\widehat{SE}$ & NA & $0.0019$ & $0.0017$ & $0.0015$ & $0.0014$ \\ 
		& CP & NA & $0.925$ & $0.941$ & $0.945$ & $0.836$\\ \midrule 
		\multirow{4}{*}{$N=2^{22}$} & Bias & $-8.8\times 10^{-8}$ & $-2.4 \times 10^{-5}$ & $1.0\times 10^{-5}$ & $1.6\times 10^{-4}$ & $2.2\times 10^{-4}$ \\ 
		& SD & $0.0021$ & $0.0012$ & $0.0011$ & $9.6 \times 10^{-4}$ & $8.6\times 10^{-4}$ \\
		& $\widehat{SE}$ & NA & $0.0012$ & $0.0011$ & $9.6\times 10^{-4}$ & $8.5 \times 10^{-4}$ \\
		& CP & NA & $0.942$ & $0.947$ & $0.937$ & $0.940$\\ \midrule 
		& & $\tilde{\theta}_2$ & \multicolumn{4}{c}{$\hat{\theta}_2$}\\ \midrule
		\multirow{4}{*}{$N=2^{18}$} & Bias & $2.0\times 10^{-5}$ & $6.7\times 10^{-4}$ & $0.0014$ & $0.0028$ & $0.0048$ \\ 
		& SD & $0.0049$ & $0.0030$ & $0.0027$ & $0.0025$ & $0.0020$ \\
		& $\widehat{SE}$ & NA & $0.0030$ & $0.0027$ & $0.0024$ & $0.0020$ \\ 
		& CP & NA  &$0.930$ & $0.908$& $0.764$ & $0.322$ \\ \midrule 
		\multirow{4}{*}{$N=2^{20}$} & Bias & $6.4\times 10^{-6}$ & $1.6\times 10^{-4}$ & $3.2\times 10^{-4}$ & $6.5\times 10^{-4}$ & $0.0016$ \\ 
		& SD & $0.0032$ & $0.0020$ & $0.0017$ & $0.0015$ & $0.0014$ \\
		& $\widehat{SE}$ & NA & $0.0019$ & $0.0017$ & $0.0015$ & $0.0014$ \\ 
		& CP & NA & $0.928$ & $0.940$ & $0.936$ & $0.786$\\ \midrule 
		\multirow{4}{*}{$N=2^{22}$} & Bias & $5.1\times 10^{-6}$ & $1.0\times 10^{-5}$ & $6.4\times 10^{-5}$ & $2.5\times 10^{-4}$ & $3.5\times 10^{-4}$ \\ 
		& SD & $0.0021$ & $0.0012$ & $0.0011$ & $9.6\times 10^{-4}$ & $8.6\times 10^{-4}$ \\
		& $\widehat{SE}$ & NA & $0.0012$ & $0.0011$ & $9.6\times 10^{-4}$ & $8.5 \times 10^{-4}$ \\
		& CP & NA & $0.941$ & $0.942$ & $0.932$ & $0.932$\\ \bottomrule
	\end{tabular}
\end{center}


\subsection{Value search estimator}

Consider the model $Y_i=1+A_i (2X_i-1)+e_i$, $i=1,\cdots,N$, where $X_i\sim N(0, 1)$, $e_i\sim N(0, 0.25)$, and $\prob(A_j=1)=0.5$. Under this model assumption, the optimal treatment rule takes the form,
\begin{eqnarray*}
	d^{\tiny{opt}}(x)=I(2x>1),
\end{eqnarray*}  
and hence $\beta^\star=2$. 

We conducted eight scenarios where $S=32,~64$, and $N=2^{24}, 2^{25}, 2^{26}$ and $2^{27}$. Due to the large sample size and limited computer memory, it is extremely slow to calculate the pooled estimators. Therefore, we only estimated the aggregated estimator $\hat{\theta}_0$. Simulation results based on 1000 replications are reported in Table \ref{tableE3}. Except for the case with $S=2^6$ and $N=2^{24}$ (where the bias is relatively large), the aggregated estimates have relatively small bias; the estimated standard errors are close to the standard deviations of the estimates; and  the coverage probability is close to the nominal level.

\begin{center}
	\captionof{table}{Bias, standard deviation (SD), standard error ($\widehat{SE}$) for $\hat{\theta}_0$, and coverage probabilities of CI for $\hat{\theta}_0$.}\label{tableE3}  
	\small
	\begin{tabular}{cccc|cccc}
		\toprule 
		& & $S=2^5$  & $S=2^6$ & & & $S=2^5$ & $S=2^6$\\ \midrule
		\multirow{4}{*}{$N=2^{24}$} & Bias & $0.0031$ & $0.0087$ & \multirow{4}{*}{$N=2^{26}$} & Bias & $-0.0013$ & $2.2 \times 10^{-4}$ \\ 
		& SD & $0.0087$ & $0.0086$ & & SD & $0.0057$ & $0.0049$ \\
		& $\widehat{SE}$ & $0.0091$ & $0.0086$ & & $\widehat{SE}$ & $0.0055$ & $0.0050$ \\ 
		& CP & $0.940$  &$0.840$ & & CP & $0.929$  &$0.939$  \\ \midrule 
		\multirow{4}{*}{$N=2^{25}$} & Bias & $-1.1\times 10^{-4}$ & $0.0031$ & \multirow{4}{*}{$N=2^{27}$} & Bias & $-0.0015$ & $-0.0014$ \\ 
		& SD & $0.0068$ & $0.0063$ & & SD & $0.0043$ & $0.0041$ \\
		& $\widehat{SE}$ & $0.0070$ & $0.0064$ & & $\widehat{SE}$ & $0.0039$ & $0.0044$ \\ 
		& CP & $0.954$  &$0.921$ & & CP & $0.932$  &$0.926$  \\ \bottomrule 
	\end{tabular}
\end{center}

\section{Tail inequality for $\hat{h}^{(j)}$}\label{secconceninequality}
In this section, we establish tail inequalities for  $\hat{\theta}^{(j)}$ and $\hat{h}^{(j)}$, which are used to construct $\tilde{h}^{(j)}$, a truncated version of $\hat{h}^{(j)}$ with tail equivalence.

\begin{thm}\label{thmconcenhath}
	Under Conditions (A1)-(A5), for sufficiently large $n$, there exists some constant $C_0$, such that
	\begin{eqnarray}\label{consistency}
		\prob(\hat{\theta}^{(j)} \notin N_\delta)\le 4\exp(-C_0 n).
	\end{eqnarray}
	Moreover, for sufficiently large $n$, there exist some constants $C_1, C_2>0$ and $N_0\ge 2$, such that
	\begin{eqnarray}\label{rateofconvergence}
		\prob(||\hat{h}^{(j)}||_2\ge x|\hat{\theta}^{(j)} \in N_\delta)\le C_2 \exp(-C_1 x^3), 
	\end{eqnarray}
for any $N_0\le x\le n^{1/3}\delta$.	
\end{thm}

\begin{remark}
	\eqref{consistency} and \eqref{rateofconvergence} can be viewed as generalization of the consistency and rate of convergence results established for cube root estimators \citep[cf. Corollary 4.2 in][]{Kim1990}. The tail probability of $||\hat{h}^{(j)}||_2$ is obtained based on the subexponential tail Assumption (A3) for $m(\cdot, \theta)$. 
\end{remark}

We represent $\hat{h}^{(j)}$ as
\begin{eqnarray*}
	\hat{h}^{(j)}=\arg\max_{h\in H_{n}} M_{n,j}(h)\equiv\arg\max_{h \in H_{n}} \left\{n^{1/6} \mathbb{G}_n^{(j)} (m_h)+n^{2/3} \Mean (m_h)\right\},
\end{eqnarray*}
where $H_{n}=\{h\in \mathbb{R}^d:n^{-1/3}h+\theta_0\in \Theta\}$, $\mathbb{G}_n^{(j)}=n^{1/2} (\mathbb{P}_{n}^{(j)}-\Mean )$ and $m_h(\cdot)=m(\cdot,\theta_0+n^{-1/3}h)-m(\cdot,\theta_0)$. Similarly define
\begin{eqnarray*}
	\tilde{h}^{(j)}=\arg\max_{h\in H_n\cap H_{\delta_n}} M_{n,j}(h)=\arg\max_{h\in H_n \cap H_{\delta_n}} \left\{n^{1/6} \mathbb{G}_n^{(j)} (m_h)+n^{2/3} \Mean (m_h)\right\},
\end{eqnarray*}
where $H_{\delta_n}=\{h:||h||_2\le \delta_n\}$. By its definition, we have $||\tilde{h}^{(j)}||_2\le \delta_n$. The following Corollaries are immediate applications of Theorem \ref{thmconcenhath}. 

\begin{coro}\label{corotildehtail}
	Under Conditions (A1)-(A5), for sufficiently large $n$, there exist some constants $N_0\ge 2$, $C_4$ and $C_5$, such that
	\begin{eqnarray}\label{tailtildeh}
		\prob(||\tilde{h}^{(j)}||_2>x)\le C_5\exp(-C_4x^3), \qquad \forall x\ge N_0.
	\end{eqnarray}
\end{coro}

The proof is straightforward by noting that for any $x\le n^{1/3} \delta$,
\begin{eqnarray*}
	\prob(||\tilde{h}^{(j)}||_2>x)&\le& \prob(||\tilde{h}^{(j)}||_2>x|\hat{\theta}^{(j)}\in N_\delta)\prob(\hat{\theta}^{(j)}\in N_\delta)+\prob(\hat{\theta}^{(j)}\notin N_\delta)\\
	&\le& C_2 \exp(-C_1 x^3)+4\exp(-C_0 n)\le C_5 \exp(-C_4 x^3).
\end{eqnarray*}
\begin{remark}
	Corollary \ref{corotildehtail} suggests that $\tilde{h}^{(j)}$ has finite moments of all orders. For any $a \in \mathbb{R}^d$ and positive integer $k$, this implies that the sequence of random variables $\lvert a^T \tilde{h}^{(j)}\rvert^k$ are uniformly integrable. This result is useful in establishing the convergence for moments of $\tilde{h}^{(j)}$ (see Corollary \ref{coromomentconverge}). 
\end{remark}

\begin{coro}\label{corotailequivalent}Under Conditions (A1)-(A5), taking $\delta_n=\max(3^{1/3}, 3^{1/3}/C_1^{1/3}) \log^{1/3} n$ where $C_1$ is defined in Theorem \ref{thmconcenhath}, then $\tilde{h}^{(j)}$ and $\hat{h}^{(j)}$ are tail equivalent. If $S=o(n^3)$, $\sum_{j=1}^S \tilde{h}^{(j)}$ and $\sum_{j=1}^S \hat{h}^{(j)}$ are also tail equivalent.
\end{coro}

Tail equivalence of $\tilde{h}^{(j)}$ and $\hat{h}^{(j)}$ follows by 
\begin{eqnarray}\label{taileqtildehat}
\prob\left(\tilde{h}^{(j)}\neq \hat{h}^{(j)}\right)=\prob\left(||\hat{h}^{(j)}||_2>\delta_n\right)\le \frac{C_2}{n^3}+4\exp(-C_0n),
\end{eqnarray}	
where the last inequality is implied by Theorem \ref{thmconcenhath}. The second assertion follows by an application of Bonferroni's inequality.

Corollary \ref{corotailequivalent} proves \eqref{tildeh-hath}. 
From now on, we take $\delta_n=\max(3^{1/3}, 3^{1/3}/C_1^{1/3}) \log^{1/3} n$. 
By \eqref{taileqtildehat}, Slutsky's Theorem implies $\tilde{h}^{(j)}\stackrel{d}{\to}h_0$. Applying Skorohod's representation Theorem \citep[cf. Section 9.4 in][]{Lahiri2006}, we have that there exist random vectors $\tilde{h}^{(j)\star}\stackrel{d}{=}\tilde{h}^{(j)}$ and $h_0^\star \stackrel{d}{=}h_0$ such that $\tilde{h}^{(j)\star}\to h_0^\star$, almost surely. This together with the uniform integrability of $||\tilde{h}^{(j)}||_2^k$ gives the following Corollary.  



\begin{coro}\label{coromomentconverge}
	Under Conditions (A1)-(A5), for any $a\in \mathbb{R}^d$ and integer $k\ge 1$, we have $\Mean \{(a^T \tilde{h}^{(j)})^k\}\to \Mean \{(a^T h_0)^k\}$ as $n \rightarrow \infty$.
\end{coro}

\begin{remark}
	Taking $k=2$, it proves \eqref{variance}. Moreover, Corollary \ref{coromomentconverge} suggests a simple scheme for estimating the covariance matrix $A\equiv\Cov(h_0)$ given in \eqref{hattheta0-theta0}. For any vector $a$, by law of large numbers, we obtain
	\begin{eqnarray*}
		\frac{1}{S}\sum_{j=1}^S(a^T \tilde{h}^{(j)})^2-\Mean \{(a^T \tilde{h}^{(1)})^2\} \stackrel{a.s.}{\to }0.
	\end{eqnarray*}
	This together with tail equivalence between $\tilde{h}^{(j)}$ and $\hat{h}^{(j)}$, and $\Mean \{(a^T \tilde{h}^{(1)})^2\} \to \Mean \{(a^T h_0)^2\}$ implies that $\sum_j (a^T \hat{h}^{(j)})^2/S$ converges to $a^T A a$. 
\end{remark}

\section{Analysis of the bias}\label{secbias}
In this section, we  control the accumulated bias in the aggregated estimator as in \eqref{bias}. 
Our method is inspired by the work of \cite{Pimentel2014}, which bounds the expectation of the argmax of a stochastic process by the difference of the expected suprema of the stochastic processes with and without a linear perturbation. To illustrate our idea, we first consider a trivial case by analyzing the bias $\Mean(h_0)$.

\subsection{Stochastic process with a linear perturbation}
Recall $h_0=\arg\max_{h\in \mathbb{R}^d} Z(h)$, where $Z(h)=G(h)-1/2h^T V h$.  Under Condition (A7), the covariance function $\Omega(h_1,h_2)$ of $G$ is equal to $\{L(h_1)+L(h_2)-L(h_1-h_2)\}/2$. Symmetry of $L(\cdot)$ implies $G(\cdot)\stackrel{d}{=}G(-\cdot)$ and $Z(\cdot)\stackrel{d}{=}Z(-\cdot)$. Hence, 
\begin{eqnarray*}
	\Mean (h_0)=\Mean \{\arg\max Z(h)\}=\frac{1}{2}\left[\Mean \{\arg\max Z(h)\}+\Mean \{\arg\max Z(-h)\}\right]=0.
\end{eqnarray*}
Here we provide an alternative but superfluous proof for this trivial case. 
Define the stochastic process with a linear perturbation
\begin{eqnarray*}
	Z^{\varepsilon, a}(h)=Z(h)+\varepsilon a^T h,
\end{eqnarray*}
for any $\varepsilon\in \mathbb{R}$ and $a \in \mathbb{R}^d$. We have
\begin{eqnarray*}
	\varepsilon a^T h_0+\sup_h Z(h)\le \sup_{h} Z^{\varepsilon, a}(h).
\end{eqnarray*}
Therefore, for any $\varepsilon >0$,
\begin{eqnarray}\label{boundmeaneq1}
	a^T h_0\le \frac{1}{\varepsilon}\left(\displaystyle \sup_{h} {Z^{\varepsilon,a} (h)}- \sup_{h} Z(h)\right), 
\end{eqnarray}
and
\begin{eqnarray}\label{boundmeaneq2}
	a^T h_0\ge \frac{1}{-\varepsilon}\left(\displaystyle \sup_{h} {Z^{-\varepsilon,a} (h)}- \sup_{h} Z(h)\right).
\end{eqnarray}
It follows from \eqref{boundmeaneq1} and \eqref{boundmeaneq2} that 
\begin{eqnarray}\label{boundmeaneq3}
	|\Mean (a^T h_0)|\le \displaystyle\frac{1}{\varepsilon} \max\left(\displaystyle |\Mean \{\sup_{h} {Z^{\varepsilon,a} (h)}\}- \Mean \{\sup_{h} Z(h)\}|, |\Mean \{\sup_{h} {Z^{-\varepsilon,a} (h)}\}- \Mean \{\sup_{h} {Z(h)}|\}\right).
\end{eqnarray}
In the next Lemma we show that the right-hand side of \eqref{boundmeaneq3} is of the order $O(\varepsilon)$ for any $a \in \mathbb{R}^d$ with $||a||_2=1$. Taking $\varepsilon\to 0$, we obtain $\Mean (a^T h_0)=0$, which implies $\Mean (h_0) = 0$.

\begin{lemma}\label{lemmaguassianlinearperturb}
	Let $X(h)=B(h)-h^T W h/2$, where $B(h)$ is a mean zero process  with stationary increments and $W$ is a positive definite matrix. Assume $\Mean \{\sup_{h\in \mathbb{R}^d} X^{\varepsilon,a}(h)\}<\infty$, where $X^{\varepsilon,a}(h)=X(h)+\varepsilon a^T h$. Then, we have
	\begin{eqnarray*}
		\sup_{||a||_2=1}\left|\Mean \{\sup_{h \in \mathbb{R}^d} X^{\varepsilon,a}(h)\}-\Mean \{\sup_{h \in \mathbb{R}^d} X(h)\}\right|=O(\varepsilon^2). 
	\end{eqnarray*}
As a result, we have $\Mean \{\arg\max_h X(h)\}=0$. 
\end{lemma}

\begin{remark}
	Lemma \ref{lemmaguassianlinearperturb} can be viewed as a generalization of Theorem 4 in \cite{Pimentel2014}. Here we only require the underlying process to have stationary increments. In addition, we allow the underlying process to be indexed by multi-dimensional parameters. 
\end{remark}


The proof of Lemma \ref{lemmaguassianlinearperturb} relies on the stationary increments property of $B$, which implies 
$$
\sup_h X^{\varepsilon,a}(h)\stackrel{d}{=}\sup_h X(h)+B(\varepsilon W^{-1} a)+\varepsilon^2 a^T W^{-1} a/2.
$$ 
In the following lemma, we prove the finiteness of $\Mean \{\sup_h Z^{\varepsilon,a}(h)\}$. 

\begin{lemma}\label{lemmasupgaussianexists}
	Under Conditions (A1) - (A7), there exist some positive constants $C_6, C_7, C_8$ and $K_0$ such that 
	\begin{eqnarray*}
		\sup_{||a||_2=1,|\varepsilon|\le 1}\prob\left( \sup_h Z^{\varepsilon,a}(h)>x\right)\le C_6 \exp(-C_7 x^2)+C_8 \exp(-x),
	\end{eqnarray*}
for any $x\ge K_0$.	As a result, for any integer $m>0$, we have
	\begin{eqnarray*}
		\sup_{||a||_2=1}\sup_{|\varepsilon|\le 1}\Mean \left[\{\sup_{h} Z^{\varepsilon,a}(h)\}^m\right] <\infty.
	\end{eqnarray*}	
\end{lemma}

\begin{remark}
	Lemma \ref{lemmasupgaussianexists} shows that not only $\sup_h Z^{\varepsilon,a}(h)$ possesses finite moments of all orders, but also has a subexponential tail. This is quite surprising since the supremum is taken on $\mathbb{R}^d$. This result is due to the rescaling property of $L(\cdot)$. 
\end{remark}

\subsection{Nonasymptotic bound for the bias}
We now establish the order of $|\Mean (a^T \tilde{h}^{(j)})|$. Define the following process with a linear drift 
\begin{eqnarray*}
	M^{\varepsilon_n,a}_{n,j} (h)=M_{n,j}(h)+\varepsilon_n a^T h,
\end{eqnarray*}
for some sequence $\varepsilon_n$. By Condition (A3), it is immediate to see $\Mean \{\sup_{h\in H_{n}\cap H_{\delta_n}} M^{\varepsilon_n, a}_{n,j}(h)\}<\infty$ for any $\varepsilon_n$ and $a$. Similar to \eqref{boundmeaneq1}, \eqref{boundmeaneq2} and \eqref{boundmeaneq3}, we can show
\begin{eqnarray}\nonumber
	&|\sqrt{S} \Mean (a^T \tilde{h})|\le \displaystyle \frac{\sqrt{S}}{\varepsilon_n} \max\left(\displaystyle |\Mean \{\sup_{h\in H_{n}\cap H_{\delta_n}} {M_{n,j}^{\varepsilon_n, a} (h)}\}- \Mean \{\sup_{h\in H_{n}\cap H_{\delta_n}} {M_{n,j}(h)}\}|,\right. \\ \label{boundmeaneq5}
	&\left.\displaystyle |\Mean\{ \sup_{h\in H_{n}\cap H_{\delta_n}} {M_{n,j}^{-\varepsilon_n, a} (h)}\}- \Mean \{\sup_{h\in H_{n}\cap H_{\delta_n}} {M_{n,j}(h)}\}|\right),
\end{eqnarray}
for any positive sequence $\varepsilon_n$. 

Since $M_{n,j}(h)$ converges weakly to $Z(h)$, the expected supremum of $M_{n,j}(h)$ and $M_{n,j}^{\varepsilon_n,a}(h)$ should be close to those of $Z(h)$ and $Z^{\varepsilon,a}(h)$, respectively. Define
\begin{eqnarray}\label{deltan}
	\Delta_n=\sup_{|\varepsilon_n|\le 1, ||a||_2=1} |\Mean\{ \sup_{h\in H_{n}\cap H_{\delta_n}}{M_{n,j}^{\varepsilon_n,a}}(h)\}-\Mean \{\sup_{h\in \mathbb{R}^d} Z^{\varepsilon_n,a}(h)\}|.
\end{eqnarray}
It follows from \eqref{boundmeaneq3} that
\begin{eqnarray}\nonumber
	&\displaystyle |\sqrt{S} \Mean (a^T \tilde{h}^{(j)})|\le \displaystyle 2\frac{\sqrt{S}}{\varepsilon_n}\Delta_n+\displaystyle \frac{\sqrt{S}}{\varepsilon_n} \max\left(\displaystyle |\Mean\{ \sup_{h\in \mathbb{R}^d} {Z^{\varepsilon_n,a} (h)}\}- \Mean\{ \sup_{h\in \mathbb{R}^d} {Z(h)}\}|,\right. \\ \label{boundmeaneq4}
	&\left.\displaystyle |\Mean \{\sup_{h\in \mathbb{R}^d} {Z^{-\varepsilon_n,a} (h)}\}- \Mean\{ \sup_{h\in \mathbb{R}^d} {Z(h)}\}|\right).
\end{eqnarray}

The second term in \eqref{boundmeaneq4} is $O(\sqrt{S}\varepsilon_n)$ by Lemma \ref{lemmaguassianlinearperturb}. The first term in \eqref{boundmeaneq4} represents the approximation error of the expected supremum of Gaussian processes, whose order will be studied in the next section. If we take $\varepsilon_n=\sqrt{\Delta_n}$, the right-hand side of \eqref{boundmeaneq4} is $O(\sqrt{S\Delta_n})$. This suggests that the asymptotic normality \eqref{hattheta0-theta0} holds as long as $S=o(\Delta_n^{-1})$. Intuitively, this implies that the number of slices $S$ cannot diverge too fast, otherwise the bias will accumulate.

\subsection{Bound for the approximation error $\Delta_n$}
To establish an upper bound for $\Delta_n$,
we adopt the techniques in \cite{Cherno2013} and \cite{Cherno2014}. Specifically, they established the nonasymptotic bound for the following difference (see P.1590 in Chernozhukov et al., 2014 or Theorem 2.1 in Chernozhukov et al, 2013):
\begin{eqnarray*}
	\left|\Mean \{g(\max_{j=1}^p S_{n,j})\}-\Mean\{ g(\max_{j=1}^p T_{n,j})\}\right|,
\end{eqnarray*}
where $g$ is a smooth function with third order derivatives, $S_{n,j}=\sum_{i=1}^n X_{i,j}$ for some mean zero random vectors $X_i = (X_{i,1},\cdots,X_{i,p})^T \in \mathbb{R}^p$, and $T_{n,j}=\sum_{i=1}^n Y_{i,j}$ for some mean zero Gaussian vectors $Y_i = (Y_{i,1},\cdots,Y_{i,p})^T \in \mathbb{R}^p$ with the same covariance matrix as $X_i$.  

Here, we improve the result in two ways. First, it's not necessary to assume $X_i$ and $Y_i$ to be mean zero. The same conclusion holds as long as $\Mean X_i=\Mean Y_i=\mu_i$ for any $\mu_i < \infty$. Second, we improve the result by taking $g$ to be the identity function as needed in our application. We summarize our result in the following lemma.


\begin{lemma}\label{lemmamaximasumsofgaussian}
	Let $X_1,\cdots,X_n$ be independent random vectors in $\mathbb{R}^p$ with finite fourth absolute moments. Define $\Mean (X_i) = \mu_i=(\mu_{i1},\cdots,\mu_{ip})^T$, $Z=\max_{1\le j\le p}\sum_{i=1}^n X_{ij}$ and $\tilde{Z}=\max_{1\le j\le p} \sum Y_{ij}$, where $Y_i = (Y_{i,1},\cdots,Y_{i,p})^T$ is distributed as $N\{\mu_i, \Mean (X_i X_i^T)\}$. Then, we have for any $\beta>0$,
	\begin{eqnarray*}
		|\Mean Z-\Mean \tilde{Z}|\le 2\beta^{-1}\log p+C\beta \{B_1+\beta (B_2+B_3B_4)+\beta^2 (B_5+B_6)\},
	\end{eqnarray*} 
	where $C$ is a constant independent of $\mu_1,\cdots,\mu_n$,
	\begin{eqnarray*}
		&& B_1=\Mean \left\{\max_{1\le j,k\le p}|\sum_{i=1}^n \tilde{X}_{ij} \tilde{X}_{ik}-n\Mean (\tilde{X}_{ij} \tilde{X}_{ik})|\right\},\\
		&&B_2=\Mean \left\{\max_{1\le j,k,l\le p} |\sum_{i=1}^n \tilde{X}_{ij} \tilde{X}_{ik} \tilde{X}_{il}-n \Mean (\tilde{X}_{ij} \tilde{X}_{ik} \tilde{X}_{il})|\right\},\\
		&& B_3=\max_{j,k}\Mean \left\{\max_{1\le l\le p} |\sum_{i=1}^n \Mean (\tilde{X}_{ij} \tilde{X}_{ik}) \tilde{X}_{il}|\right\},\quad B_4=\Mean \left(\max_{1\le j\le p} \sum_{i=1}^n |\tilde{X}_{ij}|^4\right),\\
		&& B_5=n \Mean \left\{\max_{1\le j\le p} |\tilde{X}_{1j}|^4I(\max_{1\le j\le p} |\tilde{X}_{1j}|>\frac{1}{2\beta})\right\},
	\end{eqnarray*}
 $\tilde{X}_{ij}=X_{ij}-\mu_{ij}$ and $\tilde{X}_i=X_i-\mu_i$, $j=1,\cdots,p$, $i=1,\cdots,n$.
\end{lemma}

\begin{lemma}\label{lemmasupremumappro}
	Under Conditions (A1)-(A7), we have $\Delta_n=O(n^{-1/6}\log^{4/3} n)$ for sufficiently large $n$. 
\end{lemma}
This completes the proof of Theorem \ref{thmcuberootdivideandconquer}. 

Finally, we would like to point out that our method of analyzing bias is not specific to cubic-rate M-estimators. In fact, as long as 
the limiting Gaussian process of a M-estimator has stationary increments, the bias of the aggregated estimator can be similarly bounded using our method.

\section{Discussion}
In this paper, we provide a unified inference framework for aggregated M-estimators with cubic rates obtained by the divide and conquer method. Our results demonstrate that the aggregated estimators have faster convergence rate than the original M-estimators based on pooled data and achieve the asymptotic normality when the number of groups $S$ does not grow too fast with respect to $n$, the sample size of each group. 
It remains an open question whether the rate $S = o(n^{1/6}/\log^{4/3} n)$ is optimal in general, but the rate can be improved for some special cases. For example, consider the location estimator described in Section \ref{secexample}. Using the KMT approximation \citep{KMT1975} for a stochastic process indexed by one-dimensional parameter, it can be shown that the difference between the expected supremum of $\mathbb{P}_n^{(j)} (m_h)$ and the corresponding limiting Gaussian process is $O(n^{-1/3}\log n)$. However, our theorem states that this difference is of the order $O(n^{-1/6}\log^{4/3} n)$ in a general setting. 
In addition, we have not discussed on how to choose $S$ in practice. One possible way is to treat $S$ as a tuning parameter and use some cross-validation criteria to determine $S$. This is an interesting topic that needs further investigation. 

\bibliographystyle{Chicago}
\bibliography{2}

\appendix
\section{Proofs for major results}\label{secappendix}
For notational and conceptual simplicity of the proof, without loss of generality, we avoid discussing of the measurability issue, by assuming that the parameter space $\Theta$ is countable, and that $\Theta_{\delta}=N_\delta \cap \Theta$ is dense in $N_\delta$. Consider the location estimator, for example, we can take $\Theta=\mathbb{Q}$ where $\mathbb{Q}$ denotes the set of rationals. Under this assumption, the condition $\theta_0\in \Theta$ seems somewhat strong. Therefore, we  assume $\theta_0 \in \bar{\Theta}$ instead, where $\bar{\Theta}$ denotes the closure of $\Theta$, and impose an additional assumption:

\noindent (A8.) $\prob(\sup_{\theta \in \Theta} \mathbb{P}_n^{(j)} m(\cdot, \theta)\ge \mathbb{P}_n^{(j)} m(\cdot, \theta_0))=1$ for all $j=1,\dots,S$.   

Take the location estimator as an example. Let $\Theta=\mathbb{Q}$, observe that
\begin{eqnarray*}
\sup_{\theta \in \mathbb{Q}} \mathbb{P}_n^{(j)} [\theta-1,\theta+1]\le \mathbb{P}_n^{(j)} [\theta_0-1,\theta_0+1],
\end{eqnarray*}
only when some $X_i^{(j)}=\theta_0\pm 1$. However, this happens with probability $0$, since $X_i^{(j)}$ has a density function. Assumption (A8) therefore holds. We can similarly verify (A8) for the maximum score and value search estimator. For simplicity, in the proofs below, we assume $\theta_0\in \Theta$ and hence (A8) is not needed. 	

Define the empirical process
\begin{eqnarray*}
	V_n^{(j)}(\theta)=\mathbb{G}_n^{(j)} \{m(\cdot, \theta)-m(\cdot, \theta_0)\}.
\end{eqnarray*}
To prove Theorem \ref{thmconcenhath}, we need the following two Lemmas. 

\begin{lemma}\label{lemmaboundmeanvnjtheta}
	Under Conditions (A4) and (A5), we have 
	\begin{eqnarray*}
		\Mean \{||V_n^{(j)}(\theta)||_{\Theta}\}\le c_1 \sqrt{v}\omega, 
	\end{eqnarray*}
	for some constant $c_1$.
\end{lemma}

\begin{lemma}\label{lemmaboundmeanvnjthetaskn}
	Under Conditions (A4), (A5) and (A6), there exists some constant $c_3>0$, such that
	\begin{eqnarray*}
		\Mean \{||V_n^{(j)}(\theta)||_{\Theta\cap S_{k,n}}\}\le c_3 n^{-1/6} \sqrt{k}.
	\end{eqnarray*}
\end{lemma}

Definitions of $\omega$ and $v$ are given in (A3) and (A5), respectively. We define $k$ and $S_{k,n}$ in the proof of Theorem \ref{thmconcenhath}. The proofs of Lemmas 1 and 2 are given in Appendix B. 

\subsection{Proof for Theorem \ref{thmconcenhath}}
We first prove \eqref{consistency}. 
Since $V_n^{(j)}(\theta)+\sqrt{n}\Mean \{m(\cdot, \theta)-m(\cdot, \theta_0)\}=\sqrt{n}\mathbb{P}_n^{(j)} \{m(\cdot, \theta)-m(\cdot, \theta_0)\}$, under the event $\hat{\theta}^{(j)}\notin N_\delta$, we have
\begin{eqnarray}\label{proofthm31eq1}
	\sup_{\theta \in \Theta \cap N_\delta^c} \left[V_n^{(j)}(\theta)+\sqrt{n}\Mean \{m(\cdot, \theta)-m(\cdot, \theta_0)\}\right]\ge \sup_{\theta \in \Theta} \left[V_n^{(j)}(\theta)+\sqrt{n}\Mean \{m(\cdot, \theta)-m(\cdot, \theta_0)\}\right].
\end{eqnarray}
Since $\theta_0\in \Theta$, we have
\begin{eqnarray}\label{proofthm31eq2}
	\sup_{\theta \in \Theta} \left[V_n^{(j)}(\theta)+\sqrt{n}\Mean \{m(\cdot, \theta)-m(\cdot, \theta_0)\}\right]\ge 0.
\end{eqnarray}
Combining \eqref{proofthm31eq1} and \eqref{proofthm31eq2} together, we have
\begin{eqnarray}\label{proofthm31eq3}
	&&\prob(\hat{\theta}^{(j)}\notin N_\delta)\le \prob\left(\sup_{\theta \in \Theta \cap N_\delta^c} \left[V_n^{(j)}(\theta)+\sqrt{n}\Mean \{m(\cdot, \theta)-m(\cdot, \theta_0)\}\right]\ge 0\right)\\ \nonumber
	&\le & \prob\left(\sup_{\theta \in \Theta\cap N_\delta^c} V_n^{(j)}(\theta) \ge \inf_{\theta \in \Theta \cap N_\delta^c} \sqrt{n}\Mean \{m(\cdot, \theta_0)-m(\cdot, \theta)\}\right)= \prob\left(\sup_{\theta \in \Theta\cap N_\delta^c} V_n^{(j)}(\theta) \ge \sqrt{n}\eta\right),
\end{eqnarray}
where $\eta= \Mean \{m(\cdot,\theta_0)\}-\sup_{\theta \notin N_\delta}\Mean \{m(\cdot,\theta)\}$ is positive under Condition (A1). 

Under Condition (A4), we have $\omega\equiv||M(\cdot)||_{\psi_1}<\infty$. It follows from Lemma \ref{lemmaunboundedempirical} that for all $t\ge 0$,
\begin{eqnarray}\label{proofthm31eq4}
	\prob\left(||V_n^{(j)}(\theta)||_{\Theta}\ge \frac{3}{2}\Mean ||V_n^{(j)}(\theta)||_{\Theta}+t\right)\le \exp\left(-\frac{t^2}{3\sigma^2}\right)+3\exp\left(-\frac{\sqrt{n}t}{c_0\omega}\right),
\end{eqnarray}
for some constant $c_0 > 0$, and 
\begin{eqnarray}\label{proofthm31eq5}
	\sigma^2=||\Mean \{m(\cdot, \theta)-m(\cdot,\theta_0)\}^2||_{\Theta}\le 4\Mean \{M^2(\cdot)\}\le 8\omega^2.
\end{eqnarray}
The last inequality in \eqref{proofthm31eq5} follows by Lemma \ref{lemmalpnormpsip}. By Lemma \ref{lemmaboundmeanvnjtheta}, we have, for sufficiently large $n$, $\Mean \{||V_n^{(j)}(\theta)||_{\Theta}\}\le \sqrt{n}\eta/3$. Taking $t_0=\sqrt{n}\eta/2$, the event $\sup_{\theta \in \Theta\cap N_\delta^c} V_n^{(j)}(\theta) \ge \sqrt{n}\eta$ is contained in the event 
\begin{eqnarray*}
	||V_n^{(j)}(\theta)||_{\Theta}\ge \frac{3}{2}\Mean\{ ||V_n^{(j)}(\theta)||_{\Theta}\}+t_0.
\end{eqnarray*}
Hence, it follows by \eqref{proofthm31eq3}, \eqref{proofthm31eq4} and \eqref{proofthm31eq5} that
\begin{eqnarray*}
	\prob(\hat{\theta}^{(j)}\notin N_\delta)\le \exp\left(-\frac{n \eta^2}{96\omega^2}\right)+3\exp\left(-\frac{n\eta}{2c_0\omega}\right)\le 4\exp(-C_0 n),
\end{eqnarray*}
where $C_0=\min(\eta^2/96\omega^2, \eta/2c_0\omega)$. 

It remains to show \eqref{rateofconvergence}. For any  positive integer $k\le n^{-1/3}\delta$, let $S_{k,n}$ denote the shell $\{\theta:(k-1)<n^{-1/3}||\theta-\theta_0||_2\le k\}$. Let $K$ be the smallest integer such that $K\le n^{-1/3} \delta+1$. For any integer $J>0$, the event $\{||\hat{\theta}^{(j)}-\theta_0||_2> Jn^{1/3}\}$ conditional on $\hat{\theta}^{(j)} \in N_\delta$ is contained in the event $\{\hat{\theta}^{(j)} \in \cup_{J<k\le K} S_{k,n}\}$. Similar to \eqref{proofthm31eq1}, \eqref{proofthm31eq2} and \eqref{proofthm31eq3}, we have
\begin{eqnarray}\label{proofthm31eq7}
	&&\prob(||\hat{\theta}^{(j)}-\theta_0||_2\ge n^{-1/3}J|\hat{\theta}^{(j)} \in N_\delta)\le \sum_{J<k\le K} \prob(\hat{\theta}^{(j)}\in S_{k,n})\\ \nonumber
	&\le &\sum_{J<k\le K} \prob\left(\sup_{\theta \in \Theta\cap S_{k,n}} V_n^{(j)}(\theta) \ge \inf_{\theta \in \Theta \cap S_{k,n}} \sqrt{n}\Mean \{m(\cdot, \theta_0)-m(\cdot, \theta)\}\right)
\end{eqnarray}


By the second order Taylor expansion, we have 
\begin{eqnarray}\label{proofthm31eq6}
	\Mean \{m(\cdot, \theta_0)\}-\Mean \{m(\cdot, \theta)\}=\frac{1}{2} (\theta_0-\theta)^T V(\theta^\star) (\theta_0-\theta),
\end{eqnarray}
for some $\theta^\star$ joining the line segment between $\theta_0$ and $\theta$. Under Condition (A1), we have for any $\theta^\star\in N_\delta$, $\lambda_{\min}\{V(\theta^\star)\}>0$. Since $\Mean\{ m(\cdot,\theta)\}$ is twice continuously differentiable, this suggests $c_2\equiv\inf_{\theta^\star \in N_{\delta}} \lambda_{\min}\{V(\theta^\star)\}>0$. Note that $\theta\in N_\delta$ implies $\theta^\star \in N_\delta$. Together with \eqref{proofthm31eq6}, we have
\begin{eqnarray*}
	\inf_{\theta \in \Theta \cap S_{k,n}} \sqrt{n}\Mean \{m(\cdot, \theta_0)-m(\cdot, \theta)\}\ge \inf_{\theta \in S_{k,n}} \frac{\sqrt{n}c_2}{2} ||\theta_0-\theta||_2^2\ge \frac{c_2 (k-1)^2}{2n^{1/6}}. 
\end{eqnarray*}
Combining this with \eqref{proofthm31eq7}, we obtain
\begin{eqnarray}\label{proofthm31eq8}
	\prob(n^{1/3}||\hat{\theta}^{(j)}-\theta_0||_2\ge J|\hat{\theta}^{(j)} \in N_\delta)\le \sum_{J<k\le K} \prob\left(\sup_{\theta \in \Theta \cap S_{k,n}}V_n^{(j)}(\theta)\ge \frac{c_2 (k-1)^2}{2n^{1/6}}\right).
\end{eqnarray}

To bound the right-hand side of \eqref{proofthm31eq8}, we apply the same strategy. We first provide an upper bound for $\Mean\{ ||V_n^{(j)}(\theta)||_{\Theta\cap S_{k,n}}\}$, and then apply the concentration inequality for the empirical process $||V_n^{(j)}(\theta)||_{\Theta\cap S_{k,n}}$. 
Let $N_0$ denote the smallest integer greater than $1$ such that $c_2(N_0-1)^2\ge 8c_3\sqrt{N_0}$. For any $k\ge N_0$, we have $c_2(k-1)^2\ge 8c_3\sqrt{k}$. By Lemma \ref{lemmaboundmeanvnjthetaskn}, we have
\begin{eqnarray*}
	\frac{c_2(k-1)^2}{2n^{1/6}}\ge 2\Mean ||V_n^{(j)}(\theta)||_{\Theta\cap S_{k,n}}+\frac{c_2(k-1)^2}{4n^{1/6}}\ge 2\Mean ||V_n^{(j)}(\theta)||_{\Theta\cap S_{k,n}}+\frac{c_2k^2}{16n^{1/6}}.
\end{eqnarray*} 
Therefore, for all integer $k\ge N_0$, we have
\begin{eqnarray*}
	&&\prob\left(\sup_{\theta \in \Theta \cap S_{k,n}}V_n^{(j)}(\theta)\ge \frac{c_2 (k-1)^2}{2n^{1/6}}\right)\\
	&\le &\prob\left(\sup_{\theta \in \Theta \cap S_{k,n}}|V_n^{(j)}(\theta)|-2\Mean||V_n^{(j)}(\theta)||_{\Theta\cap S_{k,n}} \ge\frac{c_2 k^2}{16n^{1/6}}\right).
\end{eqnarray*}
By Lemma \ref{lemmaunboundedempirical}, the probability on the right-hand side of the inequality is bounded above by
\begin{eqnarray}\label{proofthm31eq10}
	\exp\left(-\frac{c_4 k^4}{n^{1/3}\sigma_{n,k}^2}\right)+3\exp\left(-\frac{c_5 n^{1/3} k^2}{\omega}\right),
\end{eqnarray}
where $c_4$ and $c_5$ are some constants, and $\sigma_{n,k}^2=||\Mean \{m(\cdot,\theta)-m(\cdot,\theta_0)\}^2||_{\Theta\cap S_{n,k}}$. Since $\sigma_{n,k}^2\le c_6 n^{-1/3} k$ for some constant $c_6>0$ by Condition (A2), setting $C_1=\min[c_4/c_6, c_5/\{\omega(1+\delta)\}]$, the right-hand side of \eqref{proofthm31eq10} is bounded above by
\begin{eqnarray*}
	\exp(-C_1 k^3)+3\exp(-C_1 n^{1/3} k^2 (1+\delta) )\le 4\exp(-C_1 k^3),
\end{eqnarray*}
since $n^{-1/3} k\le (\delta+1)$. Thus, the right-hand side of \eqref{proofthm31eq8} is bounded above by
\begin{eqnarray*}
	&&\sum_{k>J} 4\exp(-C_1 k^3)\le 4\sum_{k>J} \exp(-C_1(J+1)^2 k)\\
	&\le& \frac{4\exp(-C_1(J+1)^3)}{1-\exp(-C_1(J+1)^2)}\le C_2\exp\left(-C_1(J+1)^3\right),
\end{eqnarray*}
where $C_2=4/\{1-\exp(-C_1)\}$. 

From the above discussions, we have for any $x\ge N_0$, 
\begin{eqnarray*}
	&&\prob(n^{1/3}||\hat{\theta}^{(j)}-\theta_0||_2\ge x|\hat{\theta}^{(j)} \in N_\delta)\le\prob(n^{1/3}||\hat{\theta}^{(j)}-\theta_0||_2\ge [x]|\hat{\theta}^{(j)} \in N_\delta)\\
	&\le &C_2 \exp\{-C_1 ([x]+1)\}\le C_2 \exp(-C_1 x),
\end{eqnarray*}
where $[x]$ denotes the biggest integer that is smaller or equal to $x$. This completes the proof.

\subsection{Proof for Lemma \ref{lemmaguassianlinearperturb}}
Note that
\begin{eqnarray}\nonumber
&&\Mean \{\sup_{h\in \mathbb{R}^d} X^{\varepsilon,a}(h)\}=\Mean \sup_{h\in \mathbb{R}^d} \{B(h)+\varepsilon a^T h-\frac{1}{2}h^T W h\}\\ \nonumber
&=&\Mean \left[\sup_{h \in \mathbb{R}^d} \{B(h)-\frac{1}{2}(h-\varepsilon W^{-1}a)^T W (h-\varepsilon W^{-1}a)\}\right]+\frac{1}{4}\varepsilon^2 a^T W^{-1}a\\ \label{prooflemma31eq1}
&=&\Mean [\sup_{h \in \mathbb{R}^d} \{B(h+\varepsilon W^{-1}a)-\frac{1}{2}h^T W h\}]+\frac{1}{4}\varepsilon^2 a^T W^{-1}a. 
\end{eqnarray}

By Condition (A7), the process $\tilde{B}(h)\equiv B(h+\varepsilon W^{-1} a)-B(\varepsilon W^{-1} a)$ has the same distribution as $B(h)$. Therefore, \eqref{prooflemma31eq1} is equal to
\begin{eqnarray*}
	&&\Mean [\sup_{h \in \mathbb{R}^d} \{\tilde{B}(h)-\frac{1}{2}h^T W h\}]+\Mean \{B(\varepsilon W^{-1}a)\}+\frac{1}{4}\varepsilon^2 a^T W^{-1}a\\
	&=&\Mean [\sup_{h \in \mathbb{R}^d} \{B(h)-\frac{1}{2}h^T W h\}]+\frac{1}{4}\varepsilon^2 a^T W^{-1}a=\Mean \{\sup_{h\in \mathbb{R}^d} X(h)\}+\frac{1}{4}\varepsilon^2 a^T W^{-1}a.
\end{eqnarray*}
The result therefore follows by the positive definiteness of $W$.

\subsection{Proof for Lemma \ref{lemmasupremumappro}}

To prove Lemma \ref{lemmasupremumappro}, we need the following Lemma. The definitions of $B_1, \cdots, B_6$ are given latter. The proof of Lemma \ref{lemmaB1B2B3B4B5B6} is given in Appendix B. 

\begin{lemma}\label{lemmaB1B2B3B4B5B6}
	Assume $\beta=O(n^t)$ for some $0<t<1/3$. Then, for sufficiently large $n$,
	\begin{eqnarray*}
		&& B_1=O(n^{-1/3}\log^{7/6} n\sqrt{\log N})+O(n^{-5/12}\log^{3/4} N),\\
		&& B_2=O(n^{-2/3}\log^{13/6} n\sqrt{\log N})+O(n^{-3/4}\log^{3/4} N), \\
		&& B_3=O(n^{-1/6}\sqrt{\log N}), \quad B_4=O(n^{-1/3}\log^{1/3} n),\\
		&& B_5=O(n^{-5/6}\sqrt{\log N})+O(n^{-2/3}\log^{7/3} n), \quad B_6=O(1/n).
	\end{eqnarray*}
\end{lemma}

\noindent {\textit{Proof}:} We take an $\epsilon$-net $\{h_1,\dots,h_N\}$ of the metric space $(H_{\delta_n}, ||\cdot||_2)$ with $N=N(\epsilon, H_{\delta_n}, ||\cdot||_2)$. 
For notational simplicity, let $\varepsilon = \varepsilon_n$. We have
\begin{eqnarray*}
	\Delta_n &\le& |\Mean\{ \sup_{h\in H_{n}\cap H_{\delta_n}} M_{n,j}^{\varepsilon, a}(h)\}-\Mean \{\sup_{h\in H_{n} \cap H_{\delta_n}} Z^{\varepsilon, a}(h)\}| +|\Mean \{\sup_{h\in H_{n} \cap H_{\delta_n}} Z^{\varepsilon, a}(h)\}-\Mean \{\sup_{h\in \mathbb{R}^d} Z^{\varepsilon, a}(h)\}|\\
	&\le& |\Mean\{ \max_{k=1}^N M_{n,j}^{\varepsilon, a}(h_k)\}-\Mean\{ \max_{k=1}^N Z^{\varepsilon, a}(h_k)\}|+\Mean\{ \sup_{\substack{h_1,h_2\in H_n\cap H_{\delta_n}\\ ||h_1-h_2||_2\le \epsilon}} |M_{n,j}^{\varepsilon, a}(h_1)-M_{n,j}^{\varepsilon, a}(h_2)|\}\\
	&+& \Mean\{ \sup_{\substack{h_1,h_2\in H_n\cap H_{\delta_n}\\ ||h_1-h_2||_2\le \epsilon}} |Z^{\varepsilon, a}(h_1)-Z^{\varepsilon, a}(h_2)|\}+|\Mean\{ \sup_{h\in H_{n} \cap H_{\delta_n}} Z^{\varepsilon, a}(h)\}-\Mean \{\sup_{h\in H_{n}} Z^{\varepsilon, a}(h)\}|\\
	&+&|\Mean \{\sup_{h\in H_{n}} Z^{\varepsilon, a}(h)\}-\Mean\{ \sup_{h\in \mathbb{R}^d} Z^{\varepsilon, a}(h)\}|\stackrel{\Delta}{=}\eta_1+\eta_2+\eta_3+\eta_4+\eta_5.
\end{eqnarray*}
It suffices to bound $\eta_1$, $\eta_2$, $\eta_3$, $\eta_4$ and $\eta_5$. 

\medskip

\noindent \textit{Analysis of $\eta_1$: }We first approximate $M_{n,j}^{\varepsilon, a}(h)$ by a Gaussian process $Z_{n}^{\varepsilon,a}(h)$ with the same mean $\mu_n^{\varepsilon,a}(h)$ and covariance function $\Omega_n(h_1,h_2)$ as $M_{n,j}^{\varepsilon, a}(h)$, and give an upper bound for
\begin{eqnarray*}
	I_1\stackrel{\Delta}{=}|\Mean \{\max_{k=1}^N M_{n,j}^{\varepsilon, a}(h_k)\}-\Mean \{\max_{i=1}^N Z_n^{\varepsilon,a}(h_k)\}|.
\end{eqnarray*}
Next, for the Gaussian process $Z^{\varepsilon, a}(h)$ with mean $\mu^{\varepsilon,a}(\cdot)=\lim_n \mu^{\varepsilon_n,a}(\cdot)$ and covariance function $\Omega=\lim_n \Omega_n$, we give an upper bound for
\begin{eqnarray*}
	I_2\stackrel{\Delta}{=}|\Mean \{\max_{k=1}^N Z_n^{\varepsilon,a}(h_k)\}-\Mean \{\max_{k=1}^N Z^{\varepsilon, a}(h_k)\}|.
\end{eqnarray*}
By its definition, we have $\eta_1\le I_1+I_2$. 

By Lemma \ref{lemmamaximasumsofgaussian}, we have
\begin{eqnarray}\label{I1bound}
	I_1\le 2\beta^{-1}\log N+C\beta [B_1+\beta (B_2+B_3B_4)+\beta^2 (B_5+B_6)],
\end{eqnarray} 
where the constant $C$ is independent of $a$ and $\varepsilon$, and
\begin{eqnarray*}
	&& B_1=\frac{1}{n^{2/3}}\Mean \left[\max_{1\le k,l\le N}|\sum_{i=1}^n m_{h_k} (X_i^{(j)}) m_{h_l}(X_i^{(j)})-n\Mean m_{h_k} (X_1^{(j)}) m_{h_l}(X_1^{(j)})|\right],\\
	&& B_2=\frac{1}{n} \Mean \left[\max_{1\le k,l,s\le N} |\sum_{i=1}^n m_{h_k} (X_i^{(j)}) m_{h_l}(X_i^{(j)}) m_{h_s}(X_i^{(j)})-n \Mean m_{h_k} (X_1^{(j)}) m_{h_l}(X_1^{(j)}) m_{h_s}(X_1^{(j)})|\right],\\
	&& B_3=\frac{1}{n^{1/3}}\Mean \left[\max_{1\le k\le N} |\sum_{i=1}^n m_{h_k}(X_i^{(j)})|\right], B_4=\frac{1}{n^{2/3}}\max_{1\le l,k\le N}|\Mean m_{h_k} (X_1^{(j)}) m_{h_l}(X_1^{(j)})|, \\
	&& B_5=\frac{1}{n^{4/3}}\Mean \left[\max_{1\le k\le N} \sum_{i=1}^n |m_{h_k}(X_i^{(j)})|^4\right],  B_6=\frac{1}{n^{1/3}} \Mean \left[\max_{k} |m_{h_k}(X_1^{(j)})|^4I(\max_{k} |m_{h_k}(X_1^{(j)})|>\frac{n^{1/3}}{2\beta})\right].
\end{eqnarray*}
Taking $\beta=O(n^{1/6}\log ^{-1/3} n)$, it follows by \eqref{I1bound} and Lemma \ref{lemmaB1B2B3B4B5B6} that
\begin{eqnarray}\label{I1finalbound}
I_1=O(n^{-1/6}\log^{1/3} n(\log N+\log n))+O(n^{-1/4}\log^{3/2}n\log^{3/4} N).
\end{eqnarray}

Next, we bound $I_2$. A second order Taylor expansion gives
\begin{eqnarray}\label{proofI2boundeq1}
	\mu_n^{\varepsilon,a}(h)=\varepsilon a^T h-\frac{1}{2} h^T V(\theta_0+u h)h,
\end{eqnarray}
for some $0<u<n^{-1/3}$. For all $h\in H_{\delta_n}$, it follows by Condition (A6) that 
\begin{eqnarray*}
	|h^T V(\theta_0+u h)h-h^T V h|\le ||V(\theta_0+u h)-V(\theta_0)||_2 ||h||_2^2=O(n^{-1/3}\log n).
\end{eqnarray*}
This together with \eqref{proofI2boundeq1} implies
\begin{eqnarray}\label{proofI2boundeq2}
	\sup_{|\varepsilon|\le 1, ||a||_2=1} |\mu_n^{\varepsilon,a}(h)-\mu^{\varepsilon,a}(h)|=O(n^{-1/3}\log n).
\end{eqnarray}

Consider the Gaussian process $\tilde{Z}_n^{\varepsilon,a}(h)$ indexed by $h\in \mathbb{R}^d$ with mean $\mu^{\varepsilon,a}(h)$ and covariance $\Omega_n(\cdot, \cdot)$. Since $\tilde{Z}_n^{\varepsilon,a}(h)\stackrel{d}{=}Z_n^{\varepsilon,a}(h)-\mu_n^{\varepsilon,a}(h)+\mu^{\varepsilon,a}(h)$, we have
\begin{eqnarray*}
	\Mean \{\max_{k} \tilde{Z}_n^{\varepsilon,a}(h_k)\}=\Mean[ \max_k \{Z_n^{\varepsilon,a}(h_k)-\mu_n^{\varepsilon,a}(h_k)+\mu^{\varepsilon,a}(h_k)\}].
\end{eqnarray*}
Thus, by \eqref{proofI2boundeq2},
\begin{eqnarray}\label{proofI2boundeq3}
	|\Mean\{ \max_k \tilde{Z}_n^{\varepsilon,a}(h_k)\}-\Mean\{ \max_k Z_n^{\varepsilon,a}(h_k)\}|=O(n^{-1/3}\log n).
\end{eqnarray}
Define $\Omega(h,\tilde{h})=L(h)+L(\tilde{h})-L(\tilde{h}-h)$. By Condition (A7), we have that for all $1\le j,k\le N$, 
\begin{eqnarray}\label{proofI2boundeq4}
	\max_{j,k}|\Omega(h_j,h_k)-\Omega_n(h_j,h_k)|=O\left(\frac{\log^{2/3} n}{n^{1/3}}\right).
\end{eqnarray}

Note that $(\tilde{Z}_n^{\varepsilon,a}(h_1),\dots,\tilde{Z}_n^{\varepsilon,a}(h_p))^T$ and $(Z^{\varepsilon, a}(h_1),\dots,Z^{\varepsilon, a}(h_p))^T$ are Gaussian random vectors with the same mean $\mu^{\varepsilon,a}(h)$ and covariance $(\Omega_n(h_j, h_k))_{jk}$ and $(\Omega(h_j,h_k))_{jk}$, respectively. It follows by \eqref{proofI2boundeq4} and Lemma \ref{lemmacomparisongaussian} that for any $\tilde{\beta}>0$, 
\begin{eqnarray*}
	|\Mean \{\max_k \tilde{Z}_n^{\varepsilon,a}(h_k)\}-\Mean\{ \max_k Z^{\varepsilon, a}(h_k)\}|\le 2{\tilde{\beta}}^{-1}\log N+ O\left(\tilde{\beta} \frac{\log^{2/3} n}{n^{1/3}}\right).
\end{eqnarray*}
Taking $\tilde{\beta}=n^{1/6}\log^{1/6} n$, we obtain
\begin{eqnarray*}
	|\Mean \{\max_k \tilde{Z}_n^{\varepsilon,a}(h_k)\}-\Mean \{\max_k Z^{\varepsilon, a}(h_k)\}|=O(n^{-1/6}\log^{-1/6} n\log N+n^{-1/6} \log^{5/6} n).
\end{eqnarray*}
For sufficiently large $n$, this together with \eqref{proofI2boundeq3} implies
\begin{eqnarray}\label{I2finalbound}
	I_2=O(n^{-1/6}\log^{-1/6} n\log N+n^{-1/6} \log^{5/6} n).
\end{eqnarray}
Combining \eqref{I1finalbound} with \eqref{I2finalbound}, we have for sufficiently large $n$,
\begin{eqnarray}\label{eta1finalbound}
	\eta_1=O(n^{-1/6}\log^{1/3} n(\log N+\log n))+O(n^{-1/3}\log^{3/2}n\log^{3/4} N).
\end{eqnarray}

\medskip
\noindent \textit{Analysis of $\eta_2$: }We decompose $\eta_2$ as
\begin{eqnarray*}
	\eta_2\le \Mean\{ \sup_{\substack{h_1,h_2\in H_n \cap H_{\delta_n} \\||h_1-h_2||_2\le \epsilon}} |\tilde{M}_{n,j}(h_1)-\tilde{M}_{n,j}(h_2)|\}+\sup_{\substack{h_1,h_2\in H_n\cap H_{\delta_n}\\||h_1-h_2||_2\le \epsilon}}|\mu_n^{\varepsilon,a}(h_1)-\mu_n^{\varepsilon,a}(h_2)|,
\end{eqnarray*}
where $\tilde{M}_{n,j}(h)$ is the centered process $M_{n,j}^{\varepsilon, a}(h)-\mu_n^{\varepsilon,a}(h)$, and is independent of $\varepsilon$ and $a$. For any $h_1,h_2\in H_{\delta_n}$ with $||h_1-h_2||_2\le \epsilon$, a first order Taylor expansion gives
\begin{eqnarray*}
	|\Mean (m_{h_1})-\Mean (m_{h_2})|&=&n^{1/6} \left|\Mean \{m(\cdot, \theta_0+n^{-1/3}\sqrt{\log n}h_1)-m(\cdot, \theta_0+n^{-1/3}\sqrt{\log n}h_2)\}\right|\\
	&\le &n^{-1/6} \log^{1/3} n \sup_{\theta \in \Theta} \left|\left|\frac{\partial \Mean m(\cdot,\theta)}{\partial \theta}\right|\right|_2 \epsilon.
\end{eqnarray*}
Under Condition (A1), $|\varepsilon|\le 1$ and $||a||_2\le 1$, we have
\begin{eqnarray}\label{munepsilonbound}
	\sup_{\substack{h_1,h_2\in H_{\delta_n}\\||h_1-h_2||_2\le \epsilon}}|\mu_n^{\varepsilon,a}(h_1)-\mu_n^{\varepsilon,a}(h_2)|=O(n^{-1/6}\log^{1/3} n\epsilon)+|\varepsilon a^T (h_1-h_2)|=O(\epsilon).
\end{eqnarray}
By some standard symmetrization arguments and Condition (A2), we have 
\begin{eqnarray*}
	&& \Mean [\{m_{h_1}(X_1^{(j)})-\Mean (m_{h_1})-m_{h_2}(X_2^{(j)})+\Mean (m_{h_2})\}^2]\\
	&\le& \Mean [\{m_{h_1}(X_1^{(j)})- m_{h_1}(X_2^{(j)})-m_{h_2}(X_2^{(j)})+m_{h_2}(X_2^{(j)})\}]^2\\
	&\le& 4\Mean \{(m_{h_1}-m_{h_2})^2\}\le \bar{c} \epsilon,
\end{eqnarray*}
for some constant $\bar{c}>0$. Define the metric 
$$d_h^2(h_1, h_2)=\Mean [\{m_{h_1}(X_1^{(j)})-\Mean (m_{h_1})-m_{h_2}(X_2^{(j)})+\Mean (m_{h_2})\}^2].
$$ 
We have
\begin{eqnarray}\label{proofeta2boundeq11}
	\Mean\{ \sup_{\substack{h_1,h_2\in H_n\cap H_{\delta_n}\\||h_1-h_2||_2\le \epsilon}} |\tilde{M}_{n,j}(h_1)-\tilde{M}_{n,j}(h_2)|\}\le \Mean\{ \sup_{\substack{h_1,h_2\in H_n\cap H_{\delta_n}\\ d_h(h_1,h_2)\le \sqrt{\bar{c}\epsilon} }}|\tilde{M}_{n,j}(h_1)-\tilde{M}_{n,j}(h_2)|\}.
\end{eqnarray}


Moreover, by Lemma \ref{lemmaempiricalprocessvarepsilon}, it can be shown that the right-hand side of \eqref{proofeta2boundeq11} is bounded above by
\begin{eqnarray}\label{proofeta2boundeq1}
	O(\sqrt{\epsilon v\log n})+O(v\kappa^{1/4} n^{-1/4}\log^{3/4} n)+O(v\tilde{\omega} n^{-3/8} \log^{7/8} n),
\end{eqnarray}
where
\begin{eqnarray}\nonumber
	&&\kappa= \sup_{h_1,h_2\in H_{\delta_n}} \Mean (m_{h_1}-m_{h_2}-\Mean m_{h_1}+\Mean m_{h_2})^4 \le  2||\Mean (m_h-\Mean m_h)^4||_{H_{\delta_n}}\le 32 ||\Mean m_h^4||_{H_{\delta_n}}\\ \label{proofeta2boundeq2}
	&\le & 32n^{2/3} ||4\Mean M^2(X_1^{(j)}) \{m(X_1^{(j)},\theta_0+n^{-1/3} h)-m(X_1^{(j)},\theta_0)\}^2||_{H_{\delta_n}},
\end{eqnarray}
and $\tilde{\omega}=||m_h||_{\psi_1}\le 2n^{1/6} \omega$. Decomposing the right-hand side of \eqref{proofeta2boundeq2} as
\begin{eqnarray}\label{proofeta2boundeq3}
	&&128n^{2/3} \sup_{h\in H_{\delta_n}} \Mean [M^2 I(M\le \omega\log n) \{m(\cdot,\theta_0+n^{-1/3} h)-m(\cdot,\theta_0)\}^2]\\ \label{proofeta2boundeq4}
	&+&512n^{2/3}\sup_{h\in H_{\delta_n}} \Mean \{M^4 I(M> \omega\log n)\}.
\end{eqnarray}
Under Condition (A3), the term \eqref{proofeta2boundeq3} can be bounded above by 
\begin{eqnarray*}
	128n^{2/3} \omega^2 \log^2 n \sup_{h\in H_{\delta_n}}\Mean [m(\cdot,\theta_0+n^{-1/3} h)-m(\cdot,\theta_0)]^2=O(\omega^2 n^{1/3} \log^{7/3} n),
\end{eqnarray*}
and the term \eqref{proofeta2boundeq4} is $O(\omega^4 n^{-1/3})$ by Lemma \ref{lemmapsi1bound}. This implies
\begin{eqnarray*}
	\kappa=O(\omega^2 n^{1/3} \log^{7/3} n)+O(\omega^4 n^{-1/3}).
\end{eqnarray*}
Combining this together with \eqref{proofeta2boundeq1} gives
\begin{eqnarray*}
	\Mean\{ \sup_{\substack{h_1,h_2\in H_{\delta_n}\\||h_1-h_2||_2\le \epsilon}} |\tilde{M}_{n,j}(h_1)-\tilde{M}_{n,j}(h_2)|\}=O(\sqrt{\epsilon v\log n})+O(n^{-1/6} \log^{4/3} n).
\end{eqnarray*}
This together with \eqref{munepsilonbound} implies
\begin{eqnarray}\label{eta2finalbound}
	\eta_2=O(\sqrt{\epsilon v\log n})+O(n^{-1/6} \log^{4/3} n)+O(n^{-1/6}\sqrt{\log n}\epsilon). 
\end{eqnarray}

\medskip
\noindent \textit{Analysis of $\eta_3$: } Note that $G(h)=Z^{\varepsilon, a}(h)-\mu^{\varepsilon,a}(h)$ with $\mu^{\varepsilon,a}(h) = \varepsilon a^T h-\frac{1}{2}h^T Vh$. Similar to $\eta_2$, we have 
$$
\eta_3\le \Mean\{\sup_{h_1,h_2 \in H_\epsilon} |G(h_1) - G(h_2)|\}+\sup_{h_1,h_2 \in H_\epsilon}|\mu^{\varepsilon,a}(h_1) -\mu^{\varepsilon,a}(h_2)|.
$$ 
For sufficiently large $n$, we can show
\begin{eqnarray*}
\sup_{|\varepsilon|\le 1, ||a||_2\le 1}\sup_{||h_1-h_2||_2\le \epsilon}|\mu^{\varepsilon,a}(h_1) -\mu^{\varepsilon,a}(h_2)|=O(\epsilon).
\end{eqnarray*}
In addition, we have  
\begin{eqnarray*}
	\Mean [\{G(h_1)-G(h_2)\}^2]=L(h_1-h_2)\le ||h_1-h_2||_2 \sup_{\alpha \in \mathbb{R}^d, ||\alpha||_2= 1} L(\alpha),
\end{eqnarray*}
where the last inequality follows by the rescaling property of $L$ in Condition (A7). Define $\bar{C}=\sup_{\alpha \in \mathbb{R}^d, ||\alpha||_2= 1} L(\alpha)$. As $L(\cdot)$ is continuous, we have $\bar{C}<\infty$. This suggests that the process $G(h)$ is subgaussian (see Definition \ref{defisubgaussian}) with respect to the metric $e_h(h_1,h_2)=\sqrt{\bar{C}||h_1-h_2||_2}$. Thus, by Corollary 2.2.8 in \cite{van1996}, we have
\begin{eqnarray}\label{I3firstbound}
	I_3\stackrel{\Delta}{=}\Mean\{ \sup_{e_h(h_1,h_2)\le \sqrt{\bar{C}\epsilon}} |G(h_1)-G(h_2)|\}\le K \int_{0}^{\sqrt{\bar{C}\epsilon}} \sqrt{\log D(t, H_{\delta_n}, e_h)}dt,
\end{eqnarray}
where $K$ is a positive constant and $D(\varepsilon, H_{\delta_n}, d)$ stands for the packing number, i.e, the maximum number of $\epsilon$-separated points in $H_{\delta_n}$. By the relation between the packing number and covering number, we obtain
\begin{eqnarray*}
	I_3\le K \int_{0}^{\sqrt{\bar{C}\epsilon}} \sqrt{\log N(t/2, H_{\delta_n}, e_h)}dt.
\end{eqnarray*}

Let $2\bar{T}$ be the diameter of $H_{\delta_n}$ under $||\cdot||_2$ norm. For any $t>0$, it follows by the definition of $e_h$ that
\begin{eqnarray}\label{coveringnumber}
	N\left(\frac{t}{2}, H_{\delta_n}, e_h\right)=N\left(\frac{t^2}{4\bar{C}}, H_{\delta_n}, ||\cdot||_2\right)=N\left(\frac{t^2}{4\bar{C}\bar{T}}, B_2^d, ||\cdot||_2\right),
\end{eqnarray}
where $B_2^d$ is the unit $L_2$ ball in $\mathbb{R}^d$. By Lemma \ref{lemmacoveringnumberL2ball}, the last covering number in \eqref{coveringnumber} is smaller than $(1+4\bar{C}\bar{T}/t^2)^d$. Hence we have
\begin{eqnarray}\label{proofeta3boundeq1}
	I_3\le \sqrt{d} K \int_{0}^{\sqrt{\bar{C}\epsilon}} \sqrt{\log \left(1+\frac{4\bar{C}\bar{T}}{t^2}\right)}dt.
\end{eqnarray}
We decompose the right-hand side of \eqref{proofeta3boundeq1} as $I_4+I_5$, where
\begin{eqnarray*}
	I_4&=& \sqrt{d} K \int_{0}^{\sqrt{\bar{C}\epsilon}} \sqrt{\log \left(1+\frac{4\bar{C}\bar{T}}{t^2}\right)}I(t\le \sqrt{\bar{C}\bar{T}/n})dt,\\
	I_5&=&\sqrt{d} K \int_{0}^{\sqrt{\bar{C}\epsilon}} \sqrt{\log \left(1+\frac{4\bar{C}\bar{T}}{t^2}\right)}I(t> \sqrt{\bar{C}\bar{T}/n})dt.
\end{eqnarray*}
Letting $\tilde{t}=t\sqrt{\bar{C}\bar{T}/n}$, we have
\begin{eqnarray}\nonumber
	&&I_4\le K\sqrt{\frac{d\bar{C}\bar{T}}{n}}\int_{0}^1 \sqrt{\log \left(1+\frac{4n}{\tilde{t}^2}\right)}d\tilde{t}=K\sqrt{\frac{d\bar{C}\bar{T}}{n}}\int_{0}^1 \sqrt{\log n+ \log\left( \frac{1}{n}+\frac{4}{\tilde{t}^2}\right)}d\tilde{t}\\ \label{proofeta3boundeq2}
	&\le & K\sqrt{\frac{d\bar{C}\bar{T}\log n}{n}}+K\sqrt{\frac{d\bar{C}\bar{T}}{n}}\int_0^1\sqrt{\log \left(1+\frac{4}{\tilde{t}^2}\right)}d\tilde{t}=O\left(\sqrt{\frac{d\bar{T}\log n }{n}}\right),
\end{eqnarray}
where the first inequality in \eqref{proofeta3boundeq2} is due to that $\sqrt{\log a+\log b}\le \sqrt{\log a}+\sqrt{\log b}$ for all $a,b>0$. As for $I_5$, note that when $t>\sqrt{\bar{C}\bar{T}/n}$, $\log (1+4\bar{C}\bar{T}/t^2)\le \log(1+4n)\le \log(8n)$, and therefore for sufficiently large $n$,
\begin{eqnarray*}
	I_5\le \sqrt{d\bar{C}\epsilon}K \sqrt{\log (8n)}=O(\sqrt{d\epsilon\log n}). 
\end{eqnarray*}
Combining this with \eqref{proofeta3boundeq2}, we obtain
\begin{eqnarray*}
	I_3=O\left(\sqrt{\frac{d\bar{T}\log n }{n}}\right)+O(\sqrt{d\epsilon\log n}),
\end{eqnarray*}
and hence
\begin{eqnarray}\label{proofeta3finalbound}
	\eta_3=O(n^{-1/6}\sqrt{\log n}\epsilon)+O(\sqrt{d\epsilon\log n})+O\left(\sqrt{\frac{d\bar{T}\log n }{n}}\right).
\end{eqnarray}

\medskip
\noindent \textit{Analysis of $\eta_4$: }Define $\hat{h}^{\varepsilon,a}=\arg\max_{h \in H_n} Z^{\varepsilon, a}(h)$. On the set $\hat{h}^{\varepsilon,a} \in H_{\delta_n}$, we have $$\sup_{h\in H_{n}\cap H_{\delta_n}} Z^{\varepsilon, a}(h)=\sup_{h\in H_{n}} Z^{\varepsilon, a}(h),$$ and hence
\begin{eqnarray*}
	\left|\Mean \{\sup_{h\in H_{n}\cap H_{\delta_n}} Z^{\varepsilon, a}(h)\}-\Mean\{ \sup_{h\in H_{n}} Z^{\varepsilon, a}(h)\}\right|=\left|\Mean \left\{\sup_{h\in H_{n}\cap H_{\delta_n}} Z^{\varepsilon, a}(h)-\sup_{h\in H_{n}} Z^{\varepsilon, a}(h)\right\}I(\hat{h}^{\varepsilon,a} \in H_{\delta_n}^c)\right|.
\end{eqnarray*}

By Cauchy-Swartz inequality, we have  
\begin{eqnarray}\nonumber
	\eta_4^2 &\le& \Mean \left[\{\sup_{h\in H_{n}\cap H_{\delta_n}} Z^{\varepsilon, a}(h)-\sup_{h\in H_{n}} Z^{\varepsilon, a}(h)\}^2\right] \prob(\hat{h}^{\varepsilon,a} \in H_{\delta_n}^c)\\ \label{proofeta4boundeq1}
	&\le & \Mean  \left[\{\sup_{h\in H_{n}} Z^{\varepsilon, a}(h)\}^2\right] \prob(\hat{h}^{\varepsilon,a} \in H_{\delta_n}^c),
\end{eqnarray}
where the last inequality is due to $0\le \sup_{h\in H_{n}} Z^{\varepsilon, a}(h)-\sup_{h\in H_{n}\cap H_{\delta_n}} Z^{\varepsilon, a}(h)\le \sup_{h\in H_{n}} Z^{\varepsilon, a}(h)$. By Lemma \ref{lemmasupgaussianexists}, $\Mean [ \left\{\sup_{h\in H_{n}} Z^{\varepsilon, a}(h)\right\}^2]=O(1)$. Therefore it suffices to bound $\prob(\hat{h}^{\varepsilon,a} \in H_{\delta_n}^c)$. 

Since $Z^{\varepsilon, a}(0)=0$, the event $\{\hat{h}^{\varepsilon,a} \in H_{\delta_n}^c\}$ is contained in $\mathcal{A}\stackrel{\Delta}{=}\{\sup_{h\in H_{n}\cap H_{\delta_n}^c} Z^{\varepsilon, a}(h)\ge 0\}$. Define 
\begin{eqnarray*}
	\mathcal{A}_k =\left\{\sup_{\substack{h \in H_{n}\\(k-1)\bar{T}<||h||_2\le k\bar{T}}} G(h)\ge \frac{1}{2} \bar{c} (k-1)^2  \bar{T}^2 -k\bar{T}\right\},
\end{eqnarray*}
where $\bar{c}=\lambda_{\min}(V)$. It is immediate to see that
\begin{eqnarray*}
	&& \mathcal{A}\subseteq \bigcup_{k=2}^\infty \left\{ \sup_{\substack{h \in H_{n}\\ (k-1)\bar{T}<||h||_2\le k\bar{T} }}Z^{\varepsilon, a}(h)\ge 0\right\} \\
	&\subseteq&\bigcup_{k=2}^\infty \left\{ \sup_{\substack{h \in H_{n}\\ (k-1)\bar{T}<||h||_2\le k\bar{T} }}G(h)\ge \inf_{\substack{h \in H_{n}\\ (k-1)\bar{T}<||h||_2\le k\bar{T} }} \frac{1}{2}h^T V h-\varepsilon a^T h \right\}\subseteq \bigcup_{k=2}^\infty \mathcal{A}_k,
\end{eqnarray*}
since $|\varepsilon|\le 1$ and $||a||_2=1$. 

Let $H_{k,\bar{T}}$ denote the shell $(k-1)\bar{T}< ||h||_2\le k\bar{T}$. Note that $\bar{T}=c_0 \sqrt{\log n}$ for some constant $c_0$. For sufficiently large $n$ and $k\ge 2$, we have $k\bar{T}\le \bar{c}(k-1)^2\bar{T}^2/4$. Therefore
\begin{eqnarray}\label{proofboundeta4eq1}
	\prob(\mathcal{A})\le \sum_{k=2}^\infty \prob(\mathcal{A}_k)\le \sum_{k=2}^\infty \prob\left(\sup_{h\in H_{n}\cap H_{k,\bar{T}}} G(h)\ge \frac{1}{4}\bar{c}(k-1)^2 \bar{T}^2\right).
\end{eqnarray}
Under Condition (A7), the covariance function of $G$ has the following rescaling property: $\Omega(k h_1,k h_2)=k\Omega(h_1,h_2)$ for any $k>0$. This suggests $G(k \cdot)\stackrel{d}{=}\sqrt{k}G(\cdot)$ for $k>0$. Therefore,
\begin{eqnarray}\label{proofboundeta4eq2}
	\prob\left\{\sup_{h\in H_{n}\cap H_{k,\bar{T}}} G(h)\ge \frac{1}{4}\bar{c}(k-1)^2 \bar{T}^2\right\}\le \prob\left\{\sup_{h\in H_{n}/{\sqrt{k\bar{T}}}\cap B_2^d} G(h)\ge \frac{\bar{c}(k-1)^2 \bar{T}^{3/2}}{4\sqrt{k}}\right\},
\end{eqnarray}
where $H_{n}/{\sqrt{k\bar{T}}}=\{h/\sqrt{k\bar{T}}:h\in H_{n}\}$. Similar to \eqref{I3firstbound}, \eqref{coveringnumber} and \eqref{proofeta3boundeq1}, we have
\begin{eqnarray}\label{proofboundeta4eq3}
	\Mean \{\sup_{h\in B_2^d} G(h)\}=O \left(\sqrt{d} \int_0^1 \sqrt{\log (1+\frac{4}{t^2})}dt \right)=O(1).
\end{eqnarray}
Moreover, Condition (A7) suggests 
\begin{eqnarray}\label{sigmabound}
	\sigma^2=\sup_{h\in B_2^d}\Mean\{ G(h)^2\}=\sup_{h\in B_2^d} L(h)=O(1).
\end{eqnarray}
This together with \eqref{proofboundeta4eq3} implies for sufficiently large $n$, 
\begin{eqnarray*}
	\frac{\bar{c}(k-1)^2 \bar{T}^{3/2}}{4\sqrt{k}}-\Mean\{\sup_{h\in B_2^d} G(h)\}-\frac{\sigma}{\sqrt{2\pi}}\ge \frac{\bar{c}(k-1)^{3/2} \bar{T}^{3/2}}{6}-\Mean \{\sup_{h\in B_2^d} G(h)\}-\frac{\sigma}{\sqrt{2\pi}}\ge \frac{\bar{c}(k-1)^{3/2} \bar{T}^{3/2}}{12}.
\end{eqnarray*}

By Lemma \ref{lemmagaussianprocessconcentration}, the probability on the right-hand side of \eqref{proofboundeta4eq2} is thus bounded above by
\begin{eqnarray*}
	\prob\left[\sup_{h\in H_{n}/{\sqrt{k\bar{T}}}\cap B_2^d} G(h)\ge \Mean \{\sup_{h\in B_2^d} G(h)\}+\frac{\sigma}{\sqrt{2\pi}}+ \frac{\bar{c}(k-1)^{3/2} \bar{T}^{3/2}}{12}\right]\ll \exp\left(-(k-1)\bar{T}^2/c_0\right),
\end{eqnarray*}
for sufficiently large $n$. Combining this with \eqref{proofboundeta4eq1} implies 
\begin{eqnarray*}
	\prob(\mathcal{A})\le \sum_{k=2} \exp\left\{-(k-1)\bar{T}^2/c_0 \right\}\le 2\exp(-\bar{T}^2/c_0)=2/n,
\end{eqnarray*}
and hence $\eta_4=O(n^{-1/2})$. 

\medskip
\noindent \textit{Analysis of $\eta_5$: }Take $H_0$ to be a countable subset dense in $\mathbb{R}^d$, and $0 \in H_0$. Observe that
\begin{eqnarray*}
\eta_5&=&\Mean \{\sup_{h\in\mathbb{R}^d} Z^{\varepsilon,a}(h)\}-\Mean \{\sup_{h\in H_n }Z^{\varepsilon,a}(h)\}=\left(\Mean \{\sup_{h\in\mathbb{R}^d} Z^{\varepsilon,a}(h)\}-\Mean \{\sup_{h\in H_0}Z^{\varepsilon,a}(h)\}\right)\\
&+&\left(\Mean \{\sup_{h\in H_0} Z^{\varepsilon,a}(h)\}-\Mean \{\sup_{h\in H_n}Z^{\varepsilon,a}(h)\}\right)=I_6+I_7,
\end{eqnarray*}
it suffices to provide upper bounds for $I_6$ and $I_7$. We first claim for any $|\varepsilon|\le 1$, $||a||_2\le 1$, $I_6=0$.
Define $H_k=\{h:||h||_2\le k\}$ for positive integer $k$. Assume for any $k>0$, the following holds:
\begin{eqnarray}\label{proofeta5boundeq2}
\Mean \{\sup_{h\in H_0\cap H_k} Z^{\varepsilon,a}(h)\}=\Mean \{\sup_{h\in H_k} Z^{\varepsilon,a}(h)\}.
\end{eqnarray} 
Then, obviously we have for all $k>0$,
\begin{eqnarray}\label{proofeta5boundeq3}
	\Mean \{\sup_{h\in H_0\cap H_k} Z^{\varepsilon, a}(h)\}\ge \Mean \{\sup_{h\in H_k} Z^{\varepsilon, a}(h)\}.
\end{eqnarray}
Since $Z^{\varepsilon, a}(0)=0$, we have $\sup_{h\in H_0\cap H_k} Z^{\varepsilon,a}(h)\ge 0$, $\sup_{h\in H_k} Z^{\varepsilon, a}(h)\ge 0$ for all $k$. It follows by the monotone convergence theorem that
\begin{eqnarray*}
	\lim_k \Mean\{\sup_{h\in H_0 \cap H_k} Z^{\varepsilon, a}(h)\}=\Mean \{\sup_{h\in H_0} Z^{\varepsilon, a}(h)\},\quad 
	\lim_k \Mean\{\sup_{h\in H_k} Z^{\varepsilon, a}(h)\}=\Mean \{\sup_{h\in \mathbb{R}^d} Z^{\varepsilon, a}(h)\}.
\end{eqnarray*}
These together with \eqref{proofeta5boundeq3} implies
\begin{eqnarray}\label{proofeta5boundeq4}
	\Mean \{\sup_{h\in H_{0}} Z^{\varepsilon, a}(h)\}=\lim_k \Mean\{\sup_{h\in H_0 \cap H_k} Z^{\varepsilon, a}(h)\}\ge \lim_k \Mean\{\sup_{h\in H_k} Z^{\varepsilon, a}(h)\}= \Mean \{\sup_{h\in \mathbb{R}^d} Z^{\varepsilon, a}(h)\}.
\end{eqnarray}
On the other hand, we have
\begin{eqnarray}\label{proofeta5boundeq1}
\Mean \{\sup_{h\in H_{0}} Z^{\varepsilon, a}(h)\}\le \Mean \{\sup_{h\in \mathbb{R}^d} Z^{\varepsilon, a}(h)\}.
\end{eqnarray}
Assertion $I_6=0$ therefore follows by \eqref{proofeta5boundeq4} and \eqref{proofeta5boundeq1}. 

It remains to show \eqref{proofeta5boundeq2}. Recall that $e_h(h_1,h_2)=\sqrt{\bar{C}||h_1-h_2||_2}\ge \sqrt{\Mean [G(h_1)-G(h_2)]^2}$. For any $k>0$, the diameter $T_k$ of $H_k$ under the metric $\sqrt{\Mean G^2(\cdot)}$ is finite. Similar to \eqref{proofboundeta4eq3}, fo any $k>0$, we can deduce
\begin{eqnarray}\label{proofeta5boundeq5}
	\int_0^{T_k} \left(\log N(H_k, \sqrt{\Mean G^2(\cdot)}, t)\right)^{1/2}dt<\int_{0}^{T_k}\left(\log N(H_k, e_h, t)\right)^{1/2}dt<\infty.
\end{eqnarray}
By Theorem 6.1.2 in \cite{Marcus2006}, \eqref{proofeta5boundeq5} suggests there exists another version $G'$ of $G$ which is uniformly continuous on $H_k$. Since $H_{0}\cap H_k$ is dense in $H_k$, we obtain
\begin{eqnarray}\label{proofeta5eq1}
	\sup_{h\in H_{0}\cap H_k} \left(G'(h)+\varepsilon a^T h-\frac{1}{2}h^T V h\right)=\sup_{h\in H_k} \left(G'(h)+\varepsilon a^T h-\frac{1}{2}h^T V h\right),
\end{eqnarray}
and therefore $\Mean \sup_{h\in H_{n0}\cap H_k} Z^{\varepsilon, a}(h)=\Mean \sup_{h\in \tilde{H}_n\cap H_k} Z^{\varepsilon, a}(h)$. This completes the proof of \eqref{proofeta5boundeq2}.

It remains to bound $I_7$. Define $H_{n^{1/3}\delta}=\{h:||h||_2\le n^{1/3}\delta\}$. Since $\Theta_{\delta}$ is dense in $N_{\delta}$, we have $H_{n^{1/3} \delta}\cap H_n$ is dense in $H_{n^{1/3} \delta}$. Similar to \eqref{proofeta5boundeq2}, we can show
\begin{eqnarray*}
\Mean \{\sup_{h \in H_{n^{1/3}\delta}} Z^{\varepsilon, a}(h)\}=\Mean \{ \sup_{h \in H_{n^{1/3}\delta} \cap H_n} Z^{\varepsilon,a}(h) \}.
\end{eqnarray*}
Therefore, we have
\begin{eqnarray}\nonumber
I_7&\le& \Mean \{\sup_{h\in H_0} Z^{\varepsilon,a}(h)\}-\Mean \{\sup_{h\in H_n\cap H_{n^{1/3} \delta}}Z^{\varepsilon,a}(h)\}= \Mean \{\sup_{h\in H_0} Z^{\varepsilon,a}(h)\}-\Mean \{\sup_{h\in H_{n^{1/3} \delta}}Z^{\varepsilon,a}(h)\}\\ \label{proofeta5boundeq6}
&\le& \Mean \{\sup_{h\in H_0} Z^{\varepsilon,a}(h)\}-\Mean \{\sup_{h\in H_{n^{1/3} \delta} \cap H_0}Z^{\varepsilon,a}(h)\}.
\end{eqnarray}
Using similar arguments for showing $\eta_4=O(n^{-1/2})$, we have for sufficiently large $n$,
\begin{eqnarray}\label{proofI7eq2}
	\Mean\{ \sup_{h\in H_{n^{1/3} \delta}\cap H_0} Z^{\varepsilon,a}(h)\}\ge \Mean \{\sup_{h \in H_0} Z^{\varepsilon,a}(h)\}-O(n^{-1/2}).
\end{eqnarray}
Combining \eqref{proofeta5boundeq6} and \eqref{proofI7eq2} together, we obtain $I_7\le O(n^{-1/2})$. Therefore, $\eta_5=O(n^{-1/2})$. 


\medskip
So far, we have shown 
\begin{eqnarray*}
	&&\eta_1+\eta_2+\eta_3+\eta_4+\eta_5=O(n^{-1/6}\log^{1/3} n(\log n+\log N))\\
	&&+O(n^{-1/3} \log^{3/2} n\log^{3/4} N)
	+O(\sqrt{\epsilon v\log n})+O(n^{-1/6} \log^{4/3} n)\\
	&&+O(n^{-1/6}\sqrt{\log n}\epsilon)+O(\sqrt{d\epsilon\log n})+O\left(\sqrt{\frac{d\bar{T}\log n }{n}}\right)+O(n^{-1/2}).
\end{eqnarray*}
Taking $\epsilon=1/n$, it follows by Lemma \ref{lemmacoveringnumberL2ball} that the covering number $N\le (1+2n)^d\ll (4n)^d$. This implies $\Delta_n=O(n^{-1/6}\log^{4/3} n)$. The proof is therefore completed.

\subsection{Proof for Lemma \ref{lemmasupgaussianexists}}
Define $H_{k-1,k}=\{h:k-1<||h||_2\le k\}$ and $H_k=\{h:||h||_2\le k\}$. It follows by Fatou's lemma that for any $x>0$, 
\begin{eqnarray}\label{prooflemmagaussianexisteq1}
	\prob \left\{ \sup_{h\in \mathbb{R}^d} Z^{\varepsilon,a}(h)>x\right\}\le \inf \lim_k \prob \left\{\sup_{h\in H_k} Z^{\varepsilon,a}(h)>x\right\}.
\end{eqnarray}
Note that
\begin{eqnarray*}
	\sup_{h\in H_k} Z^{\varepsilon,a}(h)=\sum_{j=1}^k \sup_{h\in H_{j-1,j}} Z^{\varepsilon,a}(h).
\end{eqnarray*}
Therefore
\begin{eqnarray}\label{prooflemmagaussianexisteq2}
	\prob \left\{\sup_{h\in H_k} Z^{\varepsilon,a}(h)>x\right\}\le \sum_{j=1}^k  \prob\left\{\sup_{h\in H_{j-1,j}} Z^{\varepsilon,a}(h)>x\right\}.
\end{eqnarray}
Combining \eqref{prooflemmagaussianexisteq1} and \eqref{prooflemmagaussianexisteq2} together, it suffices to show
\begin{eqnarray}\label{prooflemmagaussianexisteq3}
	\sum_{j=1}^\infty \prob\left\{\sup_{h\in H_{j-1,j}} Z^{\varepsilon,a}(h)>x\right\}<C_6 \exp(-C_7 x^2)+C_8 \exp(-C_9 x),
\end{eqnarray}
for some constants $C_6, C_7, C_8$ and $C_9>0$. 

Define $\bar{c}=\lambda_{\min}(V)$. Since $|\varepsilon|\le 1$ and $||a||_2=1$, we have
\begin{eqnarray*}
	\sup_{h \in H_{j-1,j}}Z^{\varepsilon,a}(h)\le \sup_{h\in H_{j-1,j}}G(h)-\frac{1}{2}\bar{c}(j-1)^2+j. 
\end{eqnarray*}
By the rescaling property of $G$, the left-hand side of \eqref{prooflemmagaussianexisteq3} is bounded above by
\begin{eqnarray}\label{prooflemmagaussianexisteq5}
	\sum_{j=1}^\infty \prob\left\{\sqrt{j}\sup_{h\in B_2^d} G(h) -\frac{1}{2}\bar{c}(j-1)^2+j>x\right\},
\end{eqnarray}
where $B_2^d$ is the unit $L_2$ ball in $\mathbb{R}^d$. It follows by \eqref{proofboundeta4eq3} and \eqref{sigmabound} that $\Mean \{\sup_{h\in B_2^d} G(h)\}=O(1)$ and $\sigma^2=\sup_{h\in B_2^d} \Mean \{G^2(h)\}=O(1)$. Define $L=\Mean \{\sup_{h\in B_2^d}G(h)\}+\sigma/\sqrt{2\pi}$. There exists some integer $J_0\ge 1$ such that for $j\ge J_0$, $\bar{c}(j-1)^2\ge 4j+\sqrt{j}L$. For all $j\ge J_0$ and $x>0$, we have
\begin{eqnarray}\label{prooflemmagaussianexisteq7}
	\prob\left\{\sqrt{j}\sup_{h\in B_2^d} G(h) -\frac{1}{2}\bar{c}(j-1)^2+j>x\right\}\le \prob\left\{\sup_{h\in B_2^d} G(h)-L>\sqrt{j}+\frac{x}{\sqrt{j}}\right\}.
\end{eqnarray}

Take $B_0$ to be a countable subset dense in $B_2^d$. Similar to \eqref{proofeta5eq1}, we can show
\begin{eqnarray}\label{prooflemmagaussianexistseq10}
\prob\left(\sup_{h\in B_2^d} G(h)= \sup_{h\in B_0}G(h) \right)=1. 
\end{eqnarray}

It follows by Lemma \ref{lemmagaussianprocessconcentration}, \eqref{prooflemmagaussianexisteq7} and \eqref{prooflemmagaussianexistseq10} that
\begin{eqnarray}\nonumber
	&&\sum_{j=J_0}^\infty \prob\left\{\sqrt{j}\sup_{h\in B_2^d} G(h) -\frac{1}{2}\bar{c}(j-1)^2+j>x\right\}\le \sum_{j=J_0}^\infty \prob\left\{\sqrt{j}\sup_{h\in B_0} G(h) -\frac{1}{2}\bar{c}(j-1)^2+j>x\right\}\\ \nonumber
	&\le& \sum_{j=J_0}^\infty \prob\left\{\sup_{h\in B_0} G(h)-L>\sqrt{j}+\frac{x}{\sqrt{j}}\right\}\le \sum_{j=J_0} \frac{1}{\sqrt{2\pi}}\exp\left\{-\frac{(\sqrt{j}+x/\sqrt{j})^2}{2\sigma^2}\right\}\\ \label{prooflemmagaussianexistseq6}
	&\le &\sum_{j=J_0} \frac{1}{\sqrt{2\pi}}\exp\left(-\frac{x}{\sigma^2}\right)\exp\left(-\frac{j}{2\sigma^2}\right)=C_8 \exp(-C_9 x),
\end{eqnarray} 
for some constants $C_8, C_9>0$. On the other hand, when $j\le J_0$, define $K_0=2J_0+2\sqrt{J_0}L$. For any $x\ge K_0$, we have $x-J_0-\sqrt{J_0}L\ge x/2$. Hence,
\begin{eqnarray*}
	\prob\left\{\sqrt{j}\sup_{h\in B_2^d} G(h) -\frac{1}{2}\bar{c}(j-1)^2+j>x\right\}\le \prob\left\{\sup_{h\in B_2^d} G(h)-L>\frac{x}{2\sqrt{J_0}}\right\},
\end{eqnarray*}
for $j\le J_0$ and $x\ge K_0$. Therefore, it follows by Lemma \ref{lemmagaussianprocessconcentration} and \eqref{prooflemmagaussianexistseq10} that
\begin{eqnarray}\label{prooflemmagaussianexisteq8}
	\sum_{j=1}^{J_0} \prob\left\{\sqrt{j}\sup_{h\in B_2^d} G(h) -\frac{1}{2}\bar{c}(j-1)^2+j>x\right\}\le \frac{J_0}{\sqrt{2\pi}} \exp\left(-\frac{x^2}{8J_0\sigma^2}\right).
\end{eqnarray}
Letting $C_6=J_0/\sqrt{2\pi}$ and $C_7=1/(8J_0\sigma^2)$, the proof is completed by \eqref{prooflemmagaussianexistseq6} and \eqref{prooflemmagaussianexisteq8}.

\subsection{Proof for Lemma \ref{lemmamaximasumsofgaussian}}
Define $\bar{\mu}=\sum_{i=1}^n \mu_i/n=(\bar{\mu}_1,\dots,\bar{\mu}_p)^T$, $\tilde{X}=\sum_i X_i-n\bar{\mu}$ and $\tilde{Y}=\sum_i Y_i-n\bar{\mu}$. Moreover, define function $F_\beta(x)=\beta^{-1}\log\{\sum_j \exp(\beta x_j+\beta \bar{\mu}_j)\}: \mathbb{R}^p\to \mathbb{R}$, which approximate the maximum by a smooth function. The approximation error can be bounded by
\begin{eqnarray*}
	\max\left( |Z-F_\beta(\tilde{X})|, |\tilde{Z}-F_\beta(\tilde{Y})|\right)\le \beta^{-1} \log p.
\end{eqnarray*}
This implies
\begin{eqnarray*}
	|\Mean Z-\Mean \tilde{Z}|&\le& \Mean \left[\max\left\{ |Z-F_\beta(\tilde{X})|, |\tilde{Z}-F_\beta(\tilde{Y})|\right\}\right]+|\Mean \{F_\beta(\tilde{X})\}-\Mean\{ F_\beta(\tilde{Y})\}|\\
	&\le &2\beta^{-1}\log p+|\Mean\{ F_\beta(\tilde{X})\}-\Mean \{F_\beta(\tilde{Y})\}|.
\end{eqnarray*}
Hence it suffices to provide an upper bound for $|\Mean\{ F_\beta(\tilde{X})\}-\Mean \{F_\beta(\tilde{Y})\}|$. Here, we adopt Stein's method. For $i=1,\dots,n$, let $\tilde{X}'_i=(\tilde{X}_{i1},\dots,\tilde{X}_{ip})^T$ be an independent copy of $\tilde{X}_i \equiv X_i - \mu_i$, $i=1,\cdots,n$. In addition, let $I$ be a uniform random variable on $\{1,\dots,n\}$, independent of all other variables. Let $\tilde{X}' =\tilde{X}-\tilde{X}_I+\tilde{X}'_I$ and
\begin{eqnarray*}
	h(x)=\int_{0}^1 \frac{1}{2t} \Mean \{F_\beta(\sqrt{t}x+\sqrt{1-t}\tilde{Y})-F_\beta(\tilde{Y})\}dt.
\end{eqnarray*}
Note that $h$ has derivatives of all orders. Using similar arguments in proof of Lemma 4.2 in \cite{Cherno2014}, we can show
\begin{eqnarray*}
	|\Mean\{ F_\beta(\tilde{X})\}-\Mean \{F_\beta(\tilde{Y})\}|\le C_0\beta B_1+n \Mean (R),
\end{eqnarray*}
for some constant $C_0>0$, and
\begin{eqnarray*}
	R=h(\tilde{X}')-h(\tilde{X})-\sum_j (\tilde{X}'-\tilde{X})^T \Delta h(\tilde{X})-\frac{1}{2}\sum_{j,k} (\tilde{X}'-\tilde{X})^T \mbox{Hess} h(\tilde{X}) (\tilde{X}'-\tilde{X}),
\end{eqnarray*}
where $\Delta h$ denotes the gradient of $h$ and $\mbox{Hess}h$ denotes the Hessian matrix. It remains to bound $\Mean (R)$. By a fourth-order Taylor expansion, we obtain
\begin{eqnarray*}
	&&|\Mean (R)|\le  \left|\Mean \left\{\frac{1}{6} \sum_{i=1}^n \sum_{j,k,l=1}^p \Delta_{ij}\Delta_{ik}\Delta_{il}\partial_j \partial_k \partial_l h(\tilde{X})\right\}\right| \\
	&+&\left|\Mean \left\{\frac{1}{6} \sum_{i=1}^n \sum_{j,k,l,m=1}^p \Delta_{ij}\Delta_{ik}\Delta_{il}\Delta_{im}(1-U)^3\partial_j \partial_k \partial_l \partial_m h(\tilde{X}+U \Delta_i)\right\}\right|\stackrel{\Delta}{=}\frac{1}{6}(\zeta_1+\zeta_2),
\end{eqnarray*}
where $\Delta_{i}=(\Delta_{i1},\dots,\Delta_{ip})^T=\tilde{X}_i'-\tilde{X}_i$, $U\sim U(0, 1)$ independent of all other variables and $\partial_j$ denotes the partial derivative with respect to $x_j$.

Next, we show the following results:
\begin{eqnarray}\label{Fbetasecondderivative}
	&\displaystyle \sum_{j,k=1}^p |\partial_j \partial_k F_\beta(x)|\le 2\beta, \\ \label{Fbetathirdderivative}
	&\displaystyle \sum_{j,k,l=1}^p |\partial_j \partial_k \partial_l F_\beta(x)|\le 6\beta^2, \\ \label{Fbetafourthderivative}
	&\displaystyle \sum_{j,k,l,m=1}^p |\partial_j \partial_k \partial_l \partial_m F_\beta(x)|\le 26 \beta^3.
\end{eqnarray}
Moreover, let $U_{jklm}(x)=\sup\{|\partial_j\partial_k\partial_l\partial_m F_\beta(x+y)|:y\in \mathbb{R}^p, |y_j|\le \beta^{-1}, 1\le \forall j\le p\}$. Then there exists some constant $C_0>0$ such that
\begin{eqnarray}\label{Fbetafourthuniformderivative}
	\displaystyle \sum_{j,k,l,m=1}^p |U_{jklm}(x)|\le C_0\beta^3.
\end{eqnarray}

The proof is similar to that of Lemma 4.3 in \cite{Cherno2014}. We only show \eqref{Fbetafourthderivative} and \eqref{Fbetafourthuniformderivative} here. Define $\pi_j(x)=\exp\{\beta (x_j+\mu_j)\}/[1+\exp\{\beta (x_j+\mu_j)\}]$. We can represent $\partial_j \partial_k \partial_l \partial_m F_\beta(x)=\beta^3 q_{jklm}(x)$, where
\begin{eqnarray*}
	&&q_{jklm}(x)= \pi_j \delta_{jk}\delta_{jl}\delta_{jm}-\pi_j \pi_m \delta_{jk}\delta_{jl}-\pi_j \pi_l\delta_{jk}(\delta_{jm}+\delta_{lm})-\pi_j \pi_k (\delta_{jl}+\delta_{kl})(\delta_{jm}+\delta_{km})\\
	&+&2\pi_j\pi_k\pi_m(\delta_{jl}+\delta_{kl})+2\pi_j\pi_m\pi_l \delta_{jk}+2\pi_j\pi_k \pi_l(\delta_{lm}+\delta_{jm}+\delta_{km})-6\pi_j \pi_k\pi_m \pi_l.
\end{eqnarray*} 
Note that $\pi_j(x)\ge 0$ and $\sum_j \pi_j(x)=1$. A direct calculation shows $\sum_{j,k,l,m}|q_{jklm}(x)|\le 26$. This implies \eqref{Fbetafourthderivative}. For any $y$ with $|y_j|\le \beta^{-1}$, $1\le j\le p$, we have $\pi_j(x+y)\le e^2 \pi_j(x)$, which proves \eqref{Fbetafourthuniformderivative}.

In addition, we have
\begin{eqnarray}\nonumber
	&&\Mean \{(\tilde{X}_{ij}'-\tilde{X}_{ij})(\tilde{X}_{ik}'-\tilde{X}_{ik})(\tilde{X}_{il}'-\tilde{X}_{il})|\tilde{X}^n\}=\Mean(\tilde{X}_{ij}'\tilde{X}_{ik}' \tilde{X}_{il}')-\tilde{X}_{ij}\tilde{X}_{ik} \tilde{X}_{il}\\ \nonumber
	&-& \Mean (\tilde{X}_{ij}' \tilde{X}_{ik}') \tilde{X}_{il}-\Mean (\tilde{X}_{ij}' \tilde{X}_{il}') \tilde{X}_{ik}-\Mean (\tilde{X}_{ik}' \tilde{X}_{il}') \tilde{X}_{ij}=W_{jkl}-\tilde{X}_{ij}\tilde{X}_{ik} \tilde{X}_{il}\\ \label{prooflemmaD3eq1}
	&-&\Phi_{jk}\tilde{X}_{il}-\Phi_{jl}\tilde{X}_{ik}-\Phi_{kl}\tilde{X}_{ij},
\end{eqnarray}
where $W_{jkl}=\Mean (\tilde{X}_{ij} \tilde{X}_{ik} \tilde{X}_{il})$ and $\Phi_{jk}=\Mean (\tilde{X}_{ij} \tilde{X}_{ik})$. It follows by \eqref{prooflemmaD3eq1} that $\zeta_1$ is bounded by
\begin{eqnarray}\nonumber
	&&\zeta_1=\left|\Mean \left[\Mean \left\{\sum_i \sum_{j,k,l} \Delta_{ij} \Delta_{ik} \Delta_{il} \partial_j \partial_k \partial_l h(\tilde{X})|\tilde{X}^n\right\}\right]\right|\\ \label{prooflemmaD3eq2}
	&\le&\left|\Mean \{\sum_{j,k,l}\sum_i (\tilde{X}_{ij} \tilde{X}_{ik} \tilde{X}_{il} -W_{jkl})\partial_j \partial_k \partial_l h(\tilde{X})\}\right|+3\left|\Mean\{ \sum_{j,k,l}\sum_i \Phi_{jk} \tilde{X}_{il} \partial_j \partial_k \partial_l h(\tilde{X})\}\right|.
\end{eqnarray}
By \eqref{Fbetathirdderivative}, the first term in \eqref{prooflemmaD3eq2} is bounded by
\begin{eqnarray}\label{firsttermzeta1}
	\Mean \left\{\max_{j,k,l}\left|\sum_{i} (\tilde{X}_{ij} \tilde{X}_{ik} \tilde{X}_{il} -W_{jkl})\right|\left|\sum_{j,k,l}\partial_j \partial_k \partial_l h(\tilde{X})\right|\right\}\le 6\beta^2 B_2.
\end{eqnarray}
Similarly, we can bound the second term in \eqref{prooflemmaD3eq2} by
\begin{eqnarray}\label{secondtermzeta1}
	\max_{j,k} \Mean \left\{\max_l \left|\sum_{i=1}^n \Phi_{jk} \tilde{X}_{il}\right| \left|\sum_{j,k,l}\partial_j \partial_k \partial_l h(\tilde{X})\right|\right\}\le 6\beta^2 B_3.
\end{eqnarray}
Combining \eqref{firsttermzeta1} and \eqref{secondtermzeta1} together, we obtain
\begin{eqnarray}\label{zeta1finalbound}
	\zeta_1\le 6\beta^2 B_2+18\beta^2 B_3.
\end{eqnarray}

Define $\chi_i=I(\max_{1\le j\le p}|\Delta_{ij}|\le \beta^{-1})$ and $\chi_i^c=1-\chi_i$. Decompose $\zeta_2$ as
\begin{eqnarray}\label{prooflemmaD3eq3}
	&&\Mean \left\{\sum_{i=1}^n \sum_{j,k,l,m=1}^p \Delta_{ij}\Delta_{ik}\Delta_{il}\Delta_{im}(1-U)^3\partial_j \partial_k \partial_l \partial_m h(\tilde{X}+U \Delta_i)\chi_i\right\}\\ \label{prooflemmaD3eq4}
	&+&\Mean \left\{\sum_{i=1}^n \sum_{j,k,l,m=1}^p \Delta_{ij}\Delta_{ik}\Delta_{il}\Delta_{im}(1-U)^3\partial_j \partial_k \partial_l \partial_m h(\tilde{X}+U \Delta_i)\chi_i^c\right\}.
\end{eqnarray}
It follows by the definition of $\chi_i$ and \eqref{Fbetafourthuniformderivative} that \eqref{prooflemmaD3eq3} is smaller than
\begin{eqnarray*}
	&&\Mean \left\{\sum_{j,k,l}\max_{i}(\chi_i |\partial_j \partial_k \partial_l \partial_m h(\tilde{X}+U \Delta_i)|)\max_{j,k,l,m}\sum_i|\Delta_{ij}\Delta_{ik}\Delta_{il}\Delta_{im}|\right\}\\
	&\le &26\beta^3 \Mean \left(\max_{j,k,l,m}\sum_i|\Delta_{ij}\Delta_{ik}\Delta_{il}\Delta_{im}|\right)\le 26\beta^3 \Mean \left(\max_j \sum_i |\Delta_{ij}^4|\right)\le 316 \beta^3 \Mean \left(\max_j \sum_i |\tilde{X}_{ij}^4|\right).
\end{eqnarray*}
By \eqref{Fbetafourthderivative}, we can show that \eqref{prooflemmaD3eq4} is bounded by
\begin{eqnarray*}
	&&\sum_i \Mean\left\{ \max_{j,k,l,m}|\Delta_{ij}\Delta_{ik}\Delta_{il}\Delta_{im} \chi_i^c| \sum_{j,k,l,m}|\partial_j \partial_k \partial_l \partial_m h(\tilde{X}+U \Delta_i)|\right\}\\
	&\le& 26\beta^3 n\Mean (\chi_1^c \max_j\Delta_{1j}^4)\le 316\beta^3 n\Mean( \chi_1^c \max_j |\tilde{X}_{ij}^4|)\le 316\beta^3 nB_5,
\end{eqnarray*}
where the last inequality follows by Chebyshev's association inequalities \citep[see Lemma 4.4,][]{Cherno2014}. This implies $\zeta_2\le 316 \beta^3 (B_4+B_5)$. Together with \eqref{zeta1finalbound}, the proof is completed.

%
%
%
%
%

\section{Proofs for Lemmas in Section \ref{secappendix}}
\subsection{Proof for Lemma \ref{lemmaboundmeanvnjtheta}}
Consider the set of functions $\{m(\cdot,\theta)-m(\cdot, \theta_0):\theta\in\mathbb{R}^d\}$. Under Condition (A4) and (A5), an application of Lemma \ref{lemmaempiricalprocess} suggests
\begin{eqnarray*}
	\Mean ||V_n^{(j)}(\theta)||_{\Theta}\le \bar{c}\sqrt{\Mean M^2(X_i^{(j)})},
\end{eqnarray*}
for some constant $\bar{c}>0$. On the other hand, it follows by Lemma \ref{lemmalpnormpsip} that $\Mean M^2(X_i^{(j)})\le 2\omega^2$. Take $c_1=\sqrt{2}\bar{c}$, the assertion therefore follows.

\subsection{Proof for Lemma \ref{lemmaboundmeanvnjthetaskn}}
Consider the set of functions $\mathcal{M}_k=\{m(\cdot,\theta)-m(\cdot, \theta_0):\theta\in \Theta_0\cap S_{k,n}\}$. By Assumption (A4), the envelope function of $\mathcal{M}_k$ is bounded by $M_{kn^{-1/3}}$. Hence, it follows by Lemma \ref{lemmaempiricalprocess} and Assumptions (A4), (A5) and (A6) that
\begin{eqnarray*}
	\Mean ||V_{n}^{(j)}(\theta)||_{\Theta\cap S_{n,k}}\le c_3 \sqrt{kn^{-1/3}},
\end{eqnarray*}
for some constant $c_3>0$. This completes the proof.

\subsection{Proof for Lemma \ref{lemmaB1B2B3B4B5B6}}
We omit the superscript $(j)$ on $X_i$ for notational convenience. Using symmetrization lemma \citep[see Lemma 2.3.1 in][]{van1996}, we obtain
\begin{eqnarray}\label{proofB1boundeq1}
	B_1 \le n^{-2/3}\Mean \left[\max_{1\le k,l\le N}|\sum_{i=1}^n m_{h_k}(X_i) m_{h_l}(X_i)\varepsilon_i|\right],
\end{eqnarray}
where $\varepsilon_i$ are i.i.d Rademacher random variables. Note that for any random variable $Y$, we have $\Mean |Y|\le ||Y||_{\psi_1}\le ||Y||_{\psi_2}/\log 2$, the right-hand side of \eqref{proofB1boundeq1} is majorized by
\begin{eqnarray}\label{proofB1boundeq2}
	&& n^{-2/3} (\log 2)^{-1} \Mean ||\max_{1\le k,l\le N}|\sum_{i=1}^n m_{h_k}(X_i) m_{h_l}(X_i)\varepsilon_i|||_{\psi_2|X}\\ \nonumber
	&\le & n^{-2/3} (\log 2)^{-1}\sqrt{2\log N} \max_{k,l} \Mean ||\sum_i m_{h_k}(X_i) m_{h_l}(X_i)\varepsilon_i||_{\psi_2|X}\\ \nonumber
	&\le & n^{-2/3} (\log 2)^{-1}\sqrt{12\log N} \Mean \max_{k,l} \sqrt{\sum_i m^2_{h_k}(X_i) m^2_{h_l}(X_i)},
\end{eqnarray}
where the Orlicz norms $||\cdot||_{\psi_2|X}$ are taken over $\varepsilon_1,\dots,\varepsilon_n$ with $X_1,\dots,X_n$ fixed, the last inequality is due to Lemma 2.2.10 in \cite{van1996}. 

Using the same arguments, we can similarly show
\begin{eqnarray}\label{proofB1boundeq3}
	&&\Mean \max_{k,l} \left|\sum_i m^2_{h_k}(X_i)m_{h_l}^2(X_i)-n\Mean m^2_{h_k}(X_i)m_{h_l}^2(X_i)\right|\\ \nonumber
	&\le & (\log 2)^{-1}\sqrt{12\log N} \Mean \max_{k,l} \sqrt{\sum_{i} m^4_{h_k}(X_i) m^4_{h_l}(X_i)}\\ \nonumber
	&\le & (\log 2)^{-1}\sqrt{12\log N} \Mean \sqrt{\sum_i M^8(X_i)}=O(\sqrt{n \log N}),
\end{eqnarray}
where the last inequality follows by Lemma \ref{lemmalpnormpsip}. On the other hand, for any $k,l$, it follows by Assumption (A4) that
\begin{eqnarray}\label{proofB1boundeq6}
	&&\Mean m^2_{h_k}(X_i)m_{h_l}^2(X_i)\le \Mean M^2(X_i) m_{h_l}^2(X_i)\\ \nonumber
	&\le & \Mean M^2(X_i) m_{h_l}^2 (X_i) I(|M(X_i)|>2\omega\log n)+ \Mean M^2(X_i) m_{h_l}^2 (X_i) I(|M(X_i)|\le 2\omega \log n)\\ \nonumber
	&\le & \Mean M^4(X_i) I(|M(X_i)|>2\omega \log n)+ 4\omega^2 \log^2 n\Mean m_{h_l}^2 (X_i)= O(\frac{1}{n})+O(\delta_n \log^2 n),
\end{eqnarray}
where the last inequality follows by Lemma \ref{lemmapsi1bound} and Assumption (A3). This together with \eqref{proofB1boundeq3} implies
\begin{eqnarray}\label{proofB1boundeq4}
	\Mean \max_{k,l} \sum_i m^2_{h_k}(X_i)m_{h_l}^2(X_i)=O(\sqrt{n\log N})+O(n^{2/3}\log^{7/3} n).
\end{eqnarray}

Combining this \eqref{proofB1boundeq4} with \eqref{proofB1boundeq1} and \eqref{proofB1boundeq2} gives
\begin{eqnarray*}
	B_1=O(n^{-1/3}\log^{7/6} n\sqrt{\log N})+O(n^{-5/12}\log^{3/4} N).
\end{eqnarray*}

Similarly, we can show $B_2=O(n^{-2/3}\log^{13/6} n\sqrt{\log N})+O(n^{-3/4}\log^{3/4} N)$.

As for $B_3$, similar to \eqref{proofB1boundeq2}, we can deduce
\begin{eqnarray*}
	B_3&\le& n^{-1/3} (\log 2)^{-1}\sqrt{12\log N}\Mean \max_{k} \sqrt{\sum_i m^2_{h_k}(X_i)}\\
	&\le& n^{-1/3} (\log 2)^{-1}\sqrt{12\log N}\Mean \sqrt{\sum_i M^2(X_i)}=O(n^{-1/6}\sqrt{\log N}).
\end{eqnarray*}

Under Assumption (A3), an application of Cauchy-Swartz inequality gives
\begin{eqnarray*}
	B_4\le \sqrt{\Mean m^2_{h_k}(X_i) \Mean m^2_{h_l}(X_i)}=O(\delta_n)=O(n^{-1/3}\log^{1/3} n).
\end{eqnarray*}

Using similar arguments in \eqref{proofB1boundeq3} and \eqref{proofB1boundeq6}, we have $B_5=O(n^{-5/6}\sqrt{\log N})+O(n^{-2/3}\log^{7/3} n)$.

Under Assumption (A4), $B_6$ is bounded above by
\begin{eqnarray}\label{proofB1boundeq5}
	B_6\le n^{-1/3} \Mean M^4(X_i) I(|M(X_i)|>\frac{n^{1/3}}{2\beta}).
\end{eqnarray}

If $\beta=O(n^t)$ for some $0<t<1/3$, for sufficiently large $n$, we have $n^{1/3}/(2\beta)\gg 2\omega \log n$, and we can deduce the right-hand side of \eqref{proofB1boundeq5} is $O(1/n)$ by Lemma \ref{lemmapsi1bound}. Therefore $B_6=O(1/n)$. 

\section{Technical lemmas and definitions}
\begin{defi}\label{defiorlicznorm}
	For any random variable $Y$, define the Orlicz norm $||Y||_{\psi_p}$ as
	\begin{eqnarray*}
		||Y||_{\psi_p}\stackrel{\Delta}{=}\inf_{C>0}\left\{\Mean \exp\left(\frac{|Y|^p}{C^p}\right)\le 2\right\}.
	\end{eqnarray*}
\end{defi}

\begin{defi}\label{defisubgaussian}
	A stochastic process $X_t$ is called sub-Gaussian with respect to the semi-metric if for any $s,t\in T$, and $x>0$, we have
	\begin{eqnarray*}
		\prob(|X_s-X_t|>x)\le 2\exp\left(-\frac{x^2}{2d^2(s,t)}\right).
	\end{eqnarray*} 
\end{defi}

\begin{lemma}\label{lemmacoveringnumberL2ball}
	Given $d\ge 1$, and $\varepsilon>0$, we have
	\begin{eqnarray*}
		N(\varepsilon, B_2^d, ||\cdot||_2)\le \left(1+\frac{2}{\varepsilon}\right)^d,
	\end{eqnarray*}
	where $B_2^d$ is the unit ball in $\mathbb{R}^d$, and $||\cdot||_2$ the Euclidean metric.
\end{lemma}

\textit{Proof: }See Lemma B.3 in \cite{zhou2009}. 

\begin{lemma}\label{lemmagaussianprocessconcentration}
	Let $\{X_t: t\in T\}$ be a real-valued mean-zero Gaussian process. Assume $T$ is a countable set and $\sigma=\sup_{z\in T}(\Mean X_z^2)^{1/2}<\infty$. Then for any $t>0$, we have
	\begin{eqnarray}\label{gaussianconcentration}
		\prob\left(\sup_{z\in T}X_z>\Mean \sup_{z\in T}X_z+\frac{\sigma}{\sqrt{2\pi}}+\sigma t\right)\le  \frac{\exp(-t^2/2)}{\sqrt{2\pi} t}.
	\end{eqnarray}
\end{lemma}

\textit{Proof: }Theorem 5.4.3 and Corollary 5.4.5 in \cite{Marcus2006} implies LHS of \eqref{gaussianconcentration} is bounded by $1-\Phi(t)$ where $\Phi$ stands for the normal cdf. The assertion follows by noting that
\begin{eqnarray*}
	1-\Phi(t)\le \frac{1}{\sqrt{2\pi}} \int_{t}^{\infty} \frac{x}{t}\exp\left(-\frac{x^2}{2}\right)\le  \frac{\exp(-t^2/2)}{\sqrt{2\pi} t}.
\end{eqnarray*}

\begin{lemma}\label{lemmalpnormpsip}
	For any random variable $X$, if $||X||_{\psi_1}<\infty$, for any integer $p\ge 2$, we have
	\begin{eqnarray*}
		\Mean |X|^p \le p! ||X||_{\psi_1}^p.
	\end{eqnarray*}
\end{lemma}

\textit{Proof: }Denote $\omega=||X||_{\psi_1}$, it follows by Taylor expansion that
\begin{eqnarray*}
	2\ge \Mean \exp\left(\frac{|X|}{\omega}\right)\ge 1+\frac{\Mean |X|^p}{p!},
\end{eqnarray*}
the assertion therefore follows.

\begin{lemma}\label{lemmapsi1bound}
	Assume $||X||_{\psi_1}<\omega$, then for any integer $k\ge 0$, and $x>1$, we have
	\begin{eqnarray*}
		\Mean |X|^k I(|X|>\omega x)\le 2\sqrt{(2k)!} \omega^k \exp(-x).
	\end{eqnarray*}
\end{lemma}

\textit{Proof: }It follows by Markov's inequality that
\begin{eqnarray}\label{prooflemmapsi1boundeq1}
	\prob(|X|>\omega x)\le \Mean \frac{\exp(|X|/\omega)}{\exp(x)}\le 2\exp(-x).
\end{eqnarray}

On the other hand, it follows by Lemma \ref{lemmalpnormpsip} that
\begin{eqnarray}\label{prooflemmapsi1boundeq2}
	\Mean |X|^{2k} \le (2k)! \omega^{2k}.
\end{eqnarray}

Using Cauchy-Swartz inequality, we get $\Mean |X|^k I(|X|>\omega x)\le \sqrt{\Mean |X|^{2k} \prob(|X|>\omega x)}$, and the assertion follows by \eqref{prooflemmapsi1boundeq1} and \eqref{prooflemmapsi1boundeq2}. 

\begin{lemma}\label{lemmaunboundedempirical}
	Let $X_1, \dots, X_n$ be i.i.d random variables with values in a measurable space $(\mathcal{S}, \mathcal{B})$, and let $\mathcal{F}$ be a countable class of measurable functions $f:\mathcal{S}\to \mathbb{R}$. Assume $\Mean f(X_i)=0$ for all $f$, and $\omega=||\sup_{f\in \mathcal{F}} |f(X_i)| ||_{\psi_1}<\infty$. Then for all $0<\eta<1$ and $\delta>0$, there exists some constant $C=C(\eta,\delta)$ such that for all $t\ge 0$,
	\begin{eqnarray*}
		&&\prob \left(\sup_{f\in\mathcal{F}} |\sum_{i=1}^n f(X_i)| \ge (1+\eta) \Mean ||\sum_{i=1}^n f(X_i)||_{\mathcal{F}}+t\right)\\
		&&\le \exp\left(-\frac{t^2}{2(1+\delta)n\sigma^2}\right)+3\exp\left(-\frac{t}{C\omega}\right),
	\end{eqnarray*}  	
	where $\sigma^2=\sup_{f\in \mathcal{F}} \Mean f^2(X_i)$.
\end{lemma}

\textit{Proof: }See Theorem 4 in \cite{Adam2008}.

\begin{lemma}\label{lemmaempiricalprocess}
	Consider the same setup as in Lemma \ref{lemmaunboundedempirical}, suppose there exists some constant $K$ and $v\ge 1$ such that $\sup_Q N(\varepsilon ||F||_{Q,2}, \mathcal{F}, L_2(Q) )\le (K/\varepsilon)^{v}$, $0\le \forall \varepsilon\le 1$, where the supremum is taken over all discrete measure $Q$ such that $0<QF^2<\infty$, $L_2(Q)$ the norm on $\mathcal{F}$ defined as
	$||f||_{Q,2}=(\int |f|^2 dQ)^{1/2}$, and $F$ the envelope function $F=\sup_{f\in \mathcal{F}} |f|$. Then, there exists some constant $C>0$, such that
	\begin{eqnarray*}
		\Mean \sup_{f\in \mathcal{F}}\left|\sum_{i=1}^n f(X_i)\right|\le C\sqrt{n \Mean F^2}.
	\end{eqnarray*}
\end{lemma}

\textit{Proof: }See Theorem 2.14.1 in \cite{van1996}.

\begin{lemma}\label{lemmaempiricalprocessvarepsilon}
	Consider the same setup as in Lemma \ref{lemmaempiricalprocess}. Assume $n\ge 4$, then
	\begin{eqnarray*}
		\Mean ||\sum_i f(X_i)||_{\mathcal{F}_\epsilon}=O(\epsilon \sqrt{vn\log n})+O(v\kappa_\epsilon^{1/4} n^{1/4}\log^{3/4} n)+O(\omega v n^{1/8} \log^{7/8} n),
	\end{eqnarray*}
	where $\kappa_\epsilon=||\Mean f^4||_{\mathcal{F}_\epsilon}$, $\mathcal{F}_\epsilon=\{f_1-f_2:f_1,f_2\in\mathcal{F}, \Mean |f_1-f_2|^2\le \epsilon^2\}$.
\end{lemma}

\textit{Proof: }
Using symmetrization inequality \citep[see Lemma 2.3.1 in][]{van1996}, we have
\begin{eqnarray}\label{prooflemmad2eq1}
	\Mean ||\frac{1}{\sqrt{n}}\sum_i f(X_i)||_{\mathcal{F}_\epsilon}\le \Mean ||\frac{1}{\sqrt{n}}\sum_{i=1}^n\varepsilon_i f(X_i)||_{\mathcal{F}_\epsilon}.
\end{eqnarray}

Using similar arguments in the proof of Theorem 2.5.2 in \cite{van1996}, the right-hand side of \eqref{prooflemmad2eq1} can be further bounded from above by
\begin{eqnarray}\label{prooflemmad2eq2}
	\Mean c_0\int_{0}^{\delta_n/||F||_n} \sup_Q \sqrt{\log N(\varepsilon ||F||_{Q,2}, \mathcal{F}, L_2(Q))}d\varepsilon ||F||_n,
\end{eqnarray}
for some constant $c_0>0$, where 
\begin{eqnarray*}
	\delta_n^2= \frac{1}{n}||\sum_{i=1}^n f^2(X_i)||_\mathcal{F_\epsilon}, \qquad \mbox{and} \qquad ||F||_n^2=\sum_{i=1}^n F^2(X_i).
\end{eqnarray*}

Under the entropy assumption, the right-hand side of \eqref{prooflemmad2eq1} is majorized by
\begin{eqnarray}\label{prooflemmad4eq2}
	\sqrt{v} \Mean \int_{0}^{\delta_n/||F||_n} \sqrt{\log \left(\frac{K}{\varepsilon}\right)}d\varepsilon ||F||_n.
\end{eqnarray}

We decompose \eqref{prooflemmad4eq2} as
\begin{eqnarray*}
	\sqrt{v} \Mean \int_{0}^{1/n} \sqrt{\log \left(\frac{K}{\varepsilon}\right)}d\varepsilon ||F||_n+\sqrt{v} \Mean \int_{1/n}^{\delta_n/||F||_n} \sqrt{\log \left(\frac{K}{\varepsilon}\right)}d\varepsilon ||F||_n \stackrel{\Delta}{=}\Mean \zeta_1+\Mean \zeta_2.
\end{eqnarray*}

Consider $\zeta_1$, we have
\begin{eqnarray}\label{prooflemmad4eq3}
	\zeta_1\le \sqrt{v} \int_{0}^{1/n} \sqrt{\log \left(\frac{K}{\varepsilon}\right)}d\varepsilon ||F||_2=\frac{1}{n}\int_0^1 \sqrt{\log n+\log \left(\frac{K}{\varepsilon'}\right)}d\varepsilon' ||F||_n,
\end{eqnarray}
where the last equality follows by a change of variable $\varepsilon'=n\varepsilon$. Note that for any $a,b\ge 0$, we have $\sqrt{a+b}\le \sqrt{a}+\sqrt{b}$ and hence the right-hand side of \eqref{prooflemmad4eq3} is bounded from above by
\begin{eqnarray*}
	\frac{\sqrt{\log n}}{n}||F||_n+\frac{1}{n}\int_{0}^1 \sqrt{\log \left(\frac{K}{\varepsilon}\right)}d\varepsilon ||F||_n\le \frac{\bar{c} \log n}{n}||F||_n,
\end{eqnarray*}
for some constant $\bar{c}$, since $\log n\ge \log 3\ge 1$. This together with \eqref{prooflemmad4eq3} and Cauchy-Swartz inequality implies
\begin{eqnarray}\label{prooflemmad4zeta1}
	\Mean \zeta_1\le \frac{\bar{c} \log n}{n} \Mean ||F||_n\le\frac{\bar{c} \log n}{n} \sqrt{\Mean F^2(X_1)} =O\left(\frac{\omega \log n}{n}\right),
\end{eqnarray}
since $\Mean |F^2|\le 2||F||_{\psi_1}=2\omega^2$ by Lemma \ref{lemmapsi1bound}. Denote $a\vee b=\max(a,b)$, since $\log(K/\varepsilon)$ is decreasing as function of $\varepsilon$, $\zeta_2$ is majorized by
\begin{eqnarray}\label{prooflemmad4eq4}
	\zeta_2\le \sqrt{v} \int_{1/n}^{(\delta_n/||F||_n)\vee (1/n)} \sqrt{\log (K n)}d\varepsilon ||F||_n\le \sqrt{v\log (Kn)} \delta_n, 
\end{eqnarray}
the last inequality holds since $a\vee b-b\le a$ for $a,b>0$. It suffices to give an upper bound for $\Mean \delta_n^2$, since 
\begin{eqnarray}\label{prooflemmad4eq11}
	\Mean \zeta_2\le \sqrt{v\log (Kn)} \sqrt{\Mean \delta_n^2},
\end{eqnarray}
by \eqref{prooflemmad4eq4} and Cauchy-Swartz inequality.

We further decompose $\Mean \delta_n^2$ as
\begin{eqnarray}\label{prooflemmad4eq6}
	\Mean \delta_n^2\le \Mean ||\mathbb{P}_n f^2-\Mean f^2||_{\mathcal{F}_\epsilon}+||\Mean f^2||_{\mathcal{F}_\epsilon}.
\end{eqnarray}

The second term is smaller than or equal to $\epsilon^2$ by definition. Using similar arguments in the proof of Theorem 2.4.3 in \cite{van1996}, the first term can be bounded by
\begin{eqnarray}\label{prooflemmad4eq7}
	4\varepsilon||F||_n^2+\Mean\sqrt{1+\log N(\varepsilon ||2F||_n^2, \mathcal{F}_\epsilon, L_1(\mathbb{P}_n))} \sqrt{\frac{6}{n}}  \sqrt{||\mathbb{P}_n f^4||_{\mathcal{F}_{\epsilon}}}.
\end{eqnarray}

For any $f_1,f_2\in\mathcal{F_\epsilon}$, we have
\begin{eqnarray*}
	\mathbb{P}_n |f_1^2-f_2^2|\le \mathbb{P}_n |f_1-f_2|4F\le 4||F||_n \sqrt{ \mathbb{P}_n |f_1-f_2|^2}.
\end{eqnarray*}

The covering number $N(\varepsilon||2F||_n^2, \mathcal{F}_\epsilon^2, L_1(\mathbb{P}_n))$ is thus bounded by 
\begin{eqnarray}\label{prooflemmad4eq8}
	N(\varepsilon||2F||_n^2, \mathcal{F}_\epsilon^2, L_1(\mathbb{P}_n))\le N(\varepsilon ||F||_n, \mathcal{F}_\epsilon, L_2(\mathbb{P}_n))\le N^2\left(\frac{\varepsilon ||F||_n}{2}, \mathcal{F}, L_2(\mathbb{P}_n)\right),
\end{eqnarray}
where $\mathcal{F}_\epsilon^2=\{(f_1-f_2)^2:f_1,f_2\in\mathcal{F}_{\epsilon}\}$.

Under the entropy assumption, we can further bound the right-hand side of \eqref{prooflemmad4eq8} by $(2K/\varepsilon)^{2v}$. Take $\varepsilon=K/n$, note that $2K/\varepsilon\le 2n\le n^2$ when $n\ge 3$, and $v\ge 1$, \eqref{prooflemmad4eq7} is bounded by
\begin{eqnarray}\label{prooflemmad4eq9}
	\frac{4K}{n}\Mean ||F||_n^2+\sqrt{\frac{6+24v\log n}{n}}\Mean \sqrt{||\mathbb{P}_n f^4||_{\mathcal{F}_\epsilon}}\le \frac{8K\omega^2}{n}+\sqrt{\frac{30v\log n}{n}}\Mean \sqrt{||\mathbb{P}_n f^4||_{\mathcal{F}_\epsilon}}.
\end{eqnarray}

Similarly, we can show $\Mean ||\mathbb{P}_n f^4||_{\mathcal{F}_\epsilon}$ is bounded by
\begin{eqnarray}\label{prooflemmad4eq10}
	\kappa_\epsilon+16 \varepsilon \Mean ||F||_n ||F^3||_n+\Mean \sqrt{1+\log N(16\varepsilon ||F||_n ||F^3||_n , \mathcal{F}_\epsilon^4, L_1(\mathbb{P}_n))}\sqrt{\frac{12||\mathbb{P}_n f^8||_{\mathcal{F}}}{n}},
\end{eqnarray}
for any $\varepsilon>0$, where $\mathcal{F}_\epsilon^4=\{(f_1-f_2)^4:f_1,f_2\in\mathcal{F}_{\epsilon}\}$, $\kappa_\epsilon=||\Mean f^4||_{\mathcal{F}_\epsilon}$. Besides, we have
\begin{eqnarray*}
	N(16\varepsilon ||F||_n ||F^3||_n , \mathcal{F}_\epsilon^4, L_1(\mathbb{P}_n))\le N(4\varepsilon ||F||_n, \mathcal{F}_\epsilon, L_2(\mathbb{P}_n))\le N^2\left(2\varepsilon ||F||_n, \mathcal{F}, L_2(\mathbb{P}_n)\right),
\end{eqnarray*}
since 
\begin{eqnarray*}
	\mathbb{P}_n |f_1^4-f_2^4|\le \mathbb{P}_n |f_1-f_2| 4|F^3| \le \mathbb{P}_n ||f_1-f_2||_n ||4 F^3||_n.
\end{eqnarray*}

An application of Cauchy-Swartz inequality and Lemma \ref{lemmalpnormpsip} implies
\begin{eqnarray*}
	\Mean \sqrt{||\mathbb{P}_n f^8||_{\mathcal{F}}}\le \Mean \sqrt{\mathbb{P}_n F^8}\le \sqrt{\Mean F^8}=O(\omega^4).
\end{eqnarray*}

Take $\varepsilon=K/n$ in \eqref{prooflemmad4eq10}, $\Mean ||\mathbb{P}_n f^4||_{\mathcal{F}_\epsilon}$ is bounded by
\begin{eqnarray*}
	\kappa_\epsilon+O\left(\frac{\Mean ||F||_n ||F^3||_n}{n}\right)+O\left(\sqrt{\frac{v\log n}{n}}\omega^4\right)=\kappa_\epsilon+O\left(\sqrt{\frac{v\log n}{n}}\omega^4\right).
\end{eqnarray*}

Combining this together with \eqref{prooflemmad4eq6}, \eqref{prooflemmad4eq7} and \eqref{prooflemmad4eq9} implies
\begin{eqnarray*}
	\Mean \delta_n^2 \le \epsilon^2+O\left(\sqrt{\frac{\kappa_\epsilon v\log n}{n}}\right)+O\left(\omega^2 v^{3/4}n^{-3/4}\log^{3/4} n\right).
\end{eqnarray*}
The assertion therefore follows together with \eqref{prooflemmad4eq11}.
\begin{lemma}\label{lemmacomparisongaussian}
	Let $X=(X_1,\dots,X_p)^T$ and $Y=(Y_1,\dots,Y_p)^T$ be Gaussian random vectors with mean $\mu^X=\mu^Y=\mu$ for some $\mu\in\mathbb{R}^p$ and covariance matrices $\Sigma^X=(\sigma_{jk}^X)_{j,k}$ and $\Sigma^Y=(\sigma_{jk}^Y)_{j,k}$, respectively. Then for every $\beta>0$ and $\mu=(\mu_1,\dots,\mu_p)^T\in \mathbb{R}^p$, we have
	\begin{eqnarray*}
		\left|\Mean \max_{1\le j\le p} X_j-\Mean \max_{1\le j\le p} Y_j\right|\le 2\beta^{-1}\log p+\beta \Delta,
	\end{eqnarray*}
	where $\Delta=\max_{j,k}|\sigma_{jk}^X-\sigma_{jk}^Y|$. 
\end{lemma}

\textit{Proof: }The proof is similar to that of Theorem 1 in \cite{Cherno2015}. We note that in \cite{Cherno2015}'s paper they require $\mu=0$. However, such assumption is redundant and the above error bound holds uniformly for all $\mu\in \mathbb{R}^p$. 

%

\end{document}